\documentclass[a4paper,12pt]{article}

\usepackage[top=30truemm,bottom=25truemm,left=25truemm,right=25truemm]{geometry}
\usepackage[pdftex]{graphicx}
\renewcommand{\title}[1]{
	\vspace*{4mm}
	\begin{center}
	\textbf{\Large #1}
	\end{center}
	\smallskip}

\renewcommand{\author}[1]{
	\vspace*{0mm}
	\begin{center}
	#1
	\end{center}
	\smallskip}
\renewcommand{\abstract}[1]{
	\begin{center}
	\parbox{13cm}{\small {\sc Abstract.} #1}
	\end{center}
	\smallskip}

\renewcommand{\thanks}[1]{{
	\renewcommand{\thefootnote}{\fnsymbol{footnote}}
	\vspace*{-5mm}
	\footnote[0]{#1}}}
\newcommand{\address}[1]{\bigskip{\small\noindent #1 \par}}
\newcommand{\email}[1]{{\small\noindent\textit{Email address}: \texttt{#1} \par}}

\usepackage{titlesec} 
\titleformat*{\section}{\large\bfseries}
\titleformat*{\subsection}{\bfseries}
\titleformat*{\subsubsection}{}
\titlespacing{\section}{0pt}{*3}{*1.5}
\titlespacing{\subsection}{0pt}{*2}{*1}
\titlespacing{\subsubsection}{0pt}{*2}{*1}

 \usepackage{amscd, amsmath, amssymb, amsthm}
 \usepackage[initials,short-journals,non-sorted-cites]{amsrefs}
 \usepackage{tikz}
\usetikzlibrary{
  knots,
  hobby,
  decorations.pathreplacing,
  shapes.geometric,
  calc
}

\usepackage{graphicx}
\usepackage{color}
\usepackage{ulem}

 \theoremstyle{plain}
 \newtheorem{thm}{Theorem}[section]
 \newtheorem{prop}[thm]{Proposition}
 
 \newtheorem{lem}[thm]{Lemma}

 \theoremstyle{definition}
 \newtheorem{defi}[thm]{Definition}
 \newtheorem{example}[thm]{Example}
 \newtheorem{rem}[thm]{Remark}

\begin{document}

\tikzset{
  knot diagram/every strand/.append style={
    ultra thick,
    black 
  },
  show curve controls/.style={
    postaction=decorate,
    decoration={show path construction,
      curveto code={
        \draw [blue, dashed]
        (\tikzinputsegmentfirst)--(\tikzinputsegmentsupporta)
        node [at end, draw, solid, red, inner sep=2pt]{};
        \draw [blue, dashed]
        (\tikzinputsegmentsupportb)--(\tikzinputsegmentlast)
        node [at start, draw, solid, red, inner sep=2pt]{}
        node [at end, fill, blue, ellipse, inner sep=2pt]{}
        ;
      }
    }
  },
  show curve endpoints/.style={
    postaction=decorate,
    decoration={show path construction,
      curveto code={
        \node [fill, blue, ellipse, inner sep=2pt] at (\tikzinputsegmentlast) {}
        ;
      }
    }
  }
}
\tikzset{->-/.style 2 args={
    postaction={decorate},
    decoration={markings, mark=at position #1 with {\arrow[thick, #2]{>}}}
    },
    ->-/.default={0.5}{}
}

\tikzset{-<-/.style 2 args={
    postaction={decorate},
    decoration={markings, mark=at position #1 with {\arrow[thick, #2]{<}}}
    },
    -<-/.default={0.5}{}
}
\title{A groupoid rack and spatial surfaces}

\author{Katsunori Arai}

 \abstract{
  A spatial surface is a compact surface embedded in the $3$-sphere.
  We assume that a spatial surface is oriented 
  and that each connected component of a spatial surface is neither a disk nor without a boundary.
  A diagram of a spatial surface is a diagram of a spatial trivalent graph that is a spine of the spatial surface.
  In this paper,
  we introduce the notion of a groupoid rack, which is used for considering colorings for diagrams of spatial surfaces in order to obtain an invariant of spatial surfaces. 
  Furthermore,
  we show that a groupoid rack has a universal property on colorings for diagrams of spatial surfaces.
 }




\section{Introduction}

A \textit{spatial surface} is a compact surface embedded in the $3$-sphere $S^{3}$.
Throughout this paper, we assume that (1) a spatial surface is oriented and that 
(2) each connected component of a spatial surface is neither a disk nor without \textcolor{black}{boundaries}. 
Two spatial surfaces are \textit{equivalent} if there exists an ambient isotopy of $S^{3}$ 
which sends one to the other.

A spatial trivalent graph is a finite trivalent graph embedded in $S^{3}$.
A diagram of a spatial trivalent graph is defined as usual in knot theory.
In \cite{Matsuzaki2021}, 
S. Matsuzaki introduced a method of presenting a spatial surface by using a diagram of a spatial trivalent graph. 
Such a diagram of a spatial trivalent graph presenting a spatial surface is called a \textit{diagram} of the spatial surface. 
He also introduced Reidemeister moves for diagrams of spatial surfaces. (See Section~\ref{sect:Colorings}.)  

A \textit{rack} \cite{Fenn-Rourke1992} 
is an algebraic \textcolor{black}{system whose axioms correspond} to two of the three Reidemeister moves in knot theory.
A \textit{multiple group rack} \cite{Ishii-Matsuzaki-Murao2020} is an algebraic structure
which is used for considering colorings for diagrams of spatial surfaces.  
\textcolor{black}{For each finite multiple group rack}, 
the number of colorings using the multiple group rack is an invariant of spatial surfaces. 
In \cite{Saito-Zappala2024}, 
another algebraic \textcolor{black}{system}, which we call a \textit{heap rack}, is defined.  
It is used for considering colorings of diagrams of spatial surfaces to obtain an invariant of spatial surfaces.   

In this paper, we introduce an algebraic \textcolor{black}{system}, a \textit{groupoid rack}, 
that can be used for colorings of diagrams of spatial surfaces in order to obtain an invariant of spatial surfaces.  
For a given finite groupoid rack, 
the number of colorings is an invariant of spatial surfaces (Theorem \ref{thm coloring invariant}). 
A multiple group rack and a heap rack are regarded as groupoid racks in our sense. 
We show that a groupoid rack has a universal property on colorings of diagrams of spatial surfaces (Theorem \ref{thm universality}).

\section{A groupoid rack}

In this section, 
we introduce the notion of a groupoid rack and give some examples of groupoid racks.

A \textit{rack} \cite{Fenn-Rourke1992} is a nonempty set $R$ with a binary operation $\ast : R \times R \to R$, $(x, y) \mapsto x \ast y$\textcolor{black}{,} satisfying the following conditions.
\begin{enumerate}
	\item[(i)] For any $y \in R$, the map $S_{y} : R \to R$, defined by $x \mapsto x \ast y$, is bijective.
	\item[(ii)] For any $x, y, z \in R$, $(x \ast y) \ast z = (x \ast z) \ast (y \ast z)$.
\end{enumerate}

For any $x, y \in R$, we denote $S_{y}^{-1}(x)$ by $x \ast^{-1} y$. 
A rack $R$ with \textcolor{black}{binary} operation $\ast$ is also referred to as a rack $(R, \ast)$. 

\textcolor{black}{A category is called a \textit{groupoid} if every morphism has an inverse.}
In this paper,
we denote the composition of morphisms $f : a \to b$ and $g : b \to c$ in a category by $fg$.

\defi{\label{def:groupoindrack}}{
  A \textit{groupoid rack $X$ associated with a groupoid $\mathcal{C}$} is the set of all morphisms of the groupoid $\mathcal{C}$ equipped with a binary operation
  $\ast : X \times X \to X$ satisfying the following conditions.
  \begin{enumerate}
    \item[(i)] For any $x \in X$ and $f, g \in X$ with $\textrm{cod}(f) = \textrm{dom}(g)$,
    $x \ast (fg) = (x \ast f) \ast g$
    and $x \ast 1_{\lambda} = x$,
    where $1_{\lambda}$ is the identity morphism of the object $\lambda$.
		\item[(ii)] For any $x, y, z \in X$, $(x \ast y) \ast z = (x \ast z) \ast (y \ast z)$.
		\item[(iii)] For any $x \in X$ and $f, g \in X$ with $\textrm{cod}(f) = \textrm{dom}(g)$,
		$\textrm{cod}(f \ast x) = \textrm{dom}(g \ast x)$ and 
    $(fg) \ast x = (f \ast x) (g \ast x)$.
  \end{enumerate}  
}

A groupoid rack $X$ associated with a groupoid $\mathcal{C}$ is a rack with \textcolor{black}{binary} operation $\ast : X \times X \to X$.  
Note that the first condition of a rack follows from (i) of Definition~\ref{def:groupoindrack}. 

We may define a groupoid rack $X$ associated with a groupoid $\mathcal{C}$ to be a rack $X$ 
such that $X$ as a set is the set of all morphisms of the groupoid $\mathcal{C}$ and 
for any $x \in X$ and $f, g \in X$ with $\textrm{cod}(f) = \textrm{dom}(g)$, 
    $x \ast (fg) = (x \ast f) \ast g$ and 
    $(fg) \ast x = (f \ast x) (g \ast x)$, and 
        $x \ast 1_{\lambda} = x$ for every object $\lambda$ of $\mathcal{C}$. 

A groupoid rack associated with a groupoid is simply called a groupoid rack. 
 
\textcolor{black}{A \textit{good involution} \cite{Kamada2007,Kamada-Oshiro2010} of a rack $R$ is an involutive map $\rho : R \to R$ such that 
for any $x, y \in R$, $\rho(x \ast y) = \rho(x) \ast y$ and  $x \ast \rho(y) = x \ast^{-1} y$.
A pair $(R, \rho)$ of a rack $R$ and a good involution $\rho$ 
is called a \textit{symmetric rack}.}  

\prop{\label{Prop:Good_involution}}{
\textcolor{black}{Let $X$ be a groupoid rack associated with a groupoid $\mathcal{C}$. 
Let $\rho : X \to X$ be the  
map sending $x$ to $x^{-1}$, where $x^{-1}$ is the inverse morphism of $x$ in the groupoid $\mathcal{C}$.  
Then, $\rho$ is a good involution of $X$.}} \upshape

\begin{proof}
\textcolor{black}{First, for any $x \in X$, we have $\rho^{2}(x) = {\rho(x)}^{-1} = (x^{-1})^{-1} = x$. Thus $\rho$ is involutive.}

  \textcolor{black}{
Next, for any $x, y \in X$, by the condition (iii) in Definition~\ref{def:groupoindrack}, we obtain $$(1_{\mathrm{dom}(x)} \ast y)(x \ast y) = (1_{\mathrm{dom}(x)} x) \ast y = x \ast y.$$
By the uniqueness of the identity $1_{\mathrm{dom}(x \ast y)}$,
it follows that $1_{\mathrm{dom}(x)} \ast y = 1_{\mathrm{dom}(x \ast y)}$.
Similarly,
$1_{\mathrm{cod}(x)} \ast y = 1_{\mathrm{dom}(x \ast y)}$.
By the condition (iii), it holds that $$(x \ast y)(\rho(x) \ast y) = (x \rho(x)) \ast y = (x x^{-1}) \ast y = 1_{\mathrm{dom}}(x) \ast y = 1_{\mathrm{dom}(x \ast y)}.$$
By the uniqueness of the inverse of $x \ast y$,
$$\rho(x) \ast y = (x \ast y)^{-1} = \rho(x \ast y).$$}

\textcolor{black}{Finally, for any $x, y \in X$, using condition (i) in Definition~\ref{def:groupoindrack}, 
we obtain $$(x \ast y) \ast \rho(y) = (x \ast y) \ast y^{-1} = x \ast (y y^{-1}) = x \ast 1_{\mathrm{dom}(y)} = x.$$
Hence, we have $S_{\rho(y)} = S_{y}^{-1}$, i.e., $x \ast \rho(y) = x \ast^{-1} y$.}

\textcolor{black}{Thus $\rho$ is a good involution of $X$.}
\end{proof}

\textcolor{black}{Using the good involution as in Proposition~\ref{Prop:Good_involution}},
a groupoid rack $X$ is regarded as a symmetric rack.

We show some examples of groupoid racks.
The first two examples below show that a multiple group rack introduced in \cite{Ishii-Matsuzaki-Murao2020} and a heap rack introduced in \cite{Saito-Zappala2024} are regarded as groupoid racks in our sense.

\example[A multiple group rack, \cite{Ishii-Matsuzaki-Murao2020}]{\label{ex MGR}}{
	Let $\left\{G_{\lambda}\right\}_{\lambda \in \Lambda}$ be a family of groups and let $e_{\lambda}$ be the identity element of $G_{\lambda}$ ($\lambda \in \Lambda$).
	A \textit{multiple group rack} $X = \bigsqcup_{\lambda \in \Lambda} G_{\lambda}$ is the disjoint union of groups $G_{\lambda}$ ($\lambda \in \Lambda$)
	with a binary operation $\ast : X \times X \to X$ satisfying the following conditions.
	\begin{enumerate}
		\item[(i)] For any $x \in X$ and $y_{1}, y_{2} \in G_{\lambda}$,
		$x \ast (y_{1} y_{2}) = (x \ast y_{1}) \ast y_{2}$ and $x \ast e_{\lambda} = x$.
		\item[(ii)] For any $x, y, z \in X$, $(x \ast y) \ast z = (x \ast z) \ast (y \ast z)$.
		\item[(iii)] For any $x \in X$ and $\lambda \in \Lambda$,
		there exists $\mu \in \Lambda$ such that 
    for any $y_{1}, y_{2} \in G_{\lambda}$,
    $(y_{1} \ast x), (y_{2} \ast x) \in G_{\mu}$ and $(y_{1} y_{2}) \ast x = (y_{1} \ast x)(y_{2} \ast x)$.
	\end{enumerate}  

	Then $X$ is a groupoid rack associated with a groupoid $\mathcal{C}$ defined \textcolor{black}{as follows}.
	\begin{itemize}
		\item $\textrm{Ob}(\mathcal{C}) = \Lambda$.
		\item For any $\lambda, \mu \in \Lambda$, $\textrm{Hom}(\lambda, \mu) = \begin{cases}
      G_{\lambda} & (\lambda = \mu)\textcolor{black}{,} \\
      \textcolor{black}{\emptyset} & (\lambda \neq \mu).
    \end{cases}$
		\item The composition $G_{\lambda} \times G_{\lambda} \to G_{\lambda}$ is defined by $(x, y) \mapsto xy$.
		\item The identity morphism of $\lambda \in \Lambda$ is the identity element $e_{\lambda}$.
		\item The inverse morphism of $x \in G_{\lambda}$ is the inverse element $x^{-1}$ of $x$ in $G_\lambda$.
	\end{itemize}
}

\example[A heap rack, \cite{Saito-Zappala2024}]{\label{ex heap rack}}{
	Let $G$ be a group.
	A binary operation $\ast : G^{2} \times G^{2} \to G^{2}$ defined by $(x, y) \ast (z, w) = (xz^{-1}w, yz^{-1}w)$ is a rack operation on $G^{2}$.
	In this paper, we call the rack $G^{2}$ with the partial product defined by $(x, y)(y, z) = (x, z)$ \textcolor{black}{$(x,y,z \in G)$} a \textit{heap rack}.

	Then a heap rack is a groupoid rack associated with the groupoid $\mathcal{C}$ defined \textcolor{black}{as follows}.
	\begin{itemize}
		\item $\textrm{Ob}(\mathcal{C}) = G$.
		\item For any $x, y \in G$, $\textrm{Hom}(x, y) = \left\{(x, y)\right\}$.
		\item The composition $\mathrm{Hom}(x, y) \times \mathrm{Hom}(y, z) \to \mathrm{Hom}(x, z)$ is defined by $((x, y), (y, z)) \mapsto (x, z)$.
		\item The identity morphism of $x \in G$ is  $(x, x)$.
		\item The inverse morphism of $(x, y) \in \textrm{Hom} (x, y)$ is $(y, x)$.
	\end{itemize}
}

The next example implies that an analogous result to Example \ref{ex heap rack} holds for racks.

\example{\label{ex rack ver. of heap rack}}{
	Let $R = (R, \ast)$ be a rack.
	A binary operation $\triangleright : R^{2} \times R^{2} \to R^{2}$ defined by $(x, y) \triangleright (z, w) = ((x \ast^{-1} z) \ast w, (x \ast^{-1} z) \ast w)$ is a rack operation on $R^{2}$.

	Then the rack $R^{2} =(R^{2}, \triangleright)$ with the partial product $(x,y)(y, z) = (x, z),\ x,y,z \in R, $ is a groupoid rack associated with the groupoid $\mathcal{C}$ defined \textcolor{black}{as follows}.
	\begin{itemize}
		\item $\textrm{Ob}(\mathcal{C}) = R$.
		\item For any $x, y \in R$, $\textrm{Hom}(x, y) = \left\{(x, y)\right\}$.
		\item The composition $\textrm{Hom}(x, y) \times \textrm{Hom}(y, z) \to \mathrm{Hom}(x, z)$ is defined by $((x, y), (y, z)) \mapsto (x, z)$.
		\item The identity morphism of $x \in R$ is  $(x, x)$.
		\item The inverse morphism of $(x, y) \in \textrm{Hom}(x, y)$ is $(y, x)$.
	\end{itemize}
}

Example \ref{ex rack ver. of heap rack} 
shows a method of constructing a groupoid rack from a given rack.

\section{Diagrams of spatial surfaces}
\label{sect:Colorings}

A \textit{spatial surface} is a compact surface embedded in $S^{3} = \mathbb{R}^{3} \sqcup \left\{\infty\right\}$.
In this paper, we assume that (1) a spatial surface is oriented and that (2) each connected component of a spatial surface is neither a disk nor without \textcolor{black}{boundaries}.
Two spatial surfaces are \textit{equivalent} if there exists an ambient isotopy of $S^{3}$ 
which sends one to the other.

A \textit{spatial trivalent graph} is a finite trivalent graph embedded in $S^{3}$.
In this paper, we assume that a trivalent graph may have loops and multiple edges, and 
a spatial trivalent graph may have some \textit{$S^{1}$-components}, i.e., circles embedded in $S^{3}$.
An $S^{1}$-component of a spatial trivalent graph is regarded as an edge of the spatial trivalent graph.
A spatial trivalent graph can be presented by a diagram as usual in knot theory.
In a diagram $D$ of a spatial trivalent graph $G$, 
an \textit{$S^{1}$-component of $D$} means a sub-diagram of $D$ which presents an $S^{1}$-component of $G$.

Let $D$ be a diagram in $\mathbb{R}^{2}$ of a spatial trivalent graph.  
Consider a spatial surface $F(D)$ obtained from $D$ as shown in  Fig.~\ref{fig construct spatial surface from spatial trivalent graph diagram}. 
More precisely, take a regular neighborhood $N(D)$ of $D$ in $\mathbb{R}^{2}$ and 
replace it locally around each crossing of $D$ with two bands  in $\mathbb{R}^{3}$ as in the rightmost of the figure. 
Then we have a compact surface embedded in $\mathbb{R}^{3}$. 
Give an orientation to the surface which is induced from the orientation of $\mathbb{R}^{2}$. Considering this oriented surface to be in $S^3 = \mathbb{R}^{3} \cup \{ \infty \}$, we have a spatial surface. It is denoted by $F(D)$ and called a \textit{spatial surface obtained from $D$}.  

\begin{figure}[h]
  \centering
  \begin{tikzpicture}[line width=1.2pt, use Hobby shortcut]
    \begin{knot}[
      consider self intersections=true,
      ignore endpoint intersections=false,
      flip crossing/.list={}
    ]
    \strand(0,2)--(0,0);
    \end{knot}
    \fill[lightgray] (1.5,2)--(1.5,0)--(2,0)--(2,2);
    \draw[->](0.5,1)--(1,1);
    \draw[->](1.5,2)--(1.5,0);
    \draw[->](2,0)--(2,2);
    \draw (1.75,0.5) node{$\circlearrowleft $};
  \end{tikzpicture}
  \quad
  \begin{tikzpicture}[line width=1.2pt, use Hobby shortcut]
    \begin{knot}[
      consider self intersections=true,
      ignore endpoint intersections=false,
      flip crossing/.list={}
    ]
    \strand(1,2)--(1,0.5);
    \strand(1,0.5)--(0,0);
    \strand(1,0.5)--(2,0);
    \end{knot}
    \draw[->] (2,1)--(2.5,1);
    \fill[lightgray] (3.5,2)--(3.5,0.625)--(4,0.625)--(4,2);
    \fill[lightgray] (2.5,0)--(3.75,0.625)--(3.5,1)--(2.5,0.5);
    \fill[lightgray] (5,0)--(3.75,0.625)--(4,1)--(5,0.5);
    \fill (1,0.5) circle(0.15);
    \draw[->] (3.5,2)--(3.5,1)--(2.5,0.5);
    \draw[->] (2.5,0)--(3.75,0.625)--(5,0);
    \draw[->] (5,0.5)--(4,1)--(4,2);
    \draw (4.5,0.45) node{$\circlearrowleft$};
  \end{tikzpicture}
  \quad
  \begin{tikzpicture}[line width=1.2pt, use Hobby shortcut]
    \begin{knot}[
      consider self intersections=true,
      ignore endpoint intersections=false,
      flip crossing/.list={},
      only when rendering/.style={}
    ]
    \strand(2,2)--(0,0);
    \strand(0,2)--(2,0);
    \end{knot}
    \draw[->] (2.5,1)--(3,1);
    \fill[lightgray] (3.5,2)--(5.5,0)--(6,0)--(4,2);
    \draw[->] (3.5,2)--(5.5,0);
    \draw[->] (6,0)--(4,2);
    \fill[lightgray] (5.5,2)--(3.5,0)--(4,0)--(6,2);
    \draw[->] (5.5,2)--(3.5,0);
    \draw[->] (4,0)--(6,2);
    \draw (5.5,0.28) node{$\circlearrowleft$};
    \draw (4,0.28) node{$\circlearrowleft$};
  \end{tikzpicture}
  \caption{A construction of a spatial surface from a diagram of a spatial trivalent graph}
  {\label{fig construct spatial surface from spatial trivalent graph diagram}}
\end{figure}

Any spatial surface is equivalent to a spatial surface obtained \textcolor{black}{along} this way, cf.~\cite{Ishii-Matsuzaki-Murao2020,Matsuzaki2021}. 

A \textit{diagram} of a spatial surface $F$ means a diagram $D$ of a spatial trivalent graph such that $F(D)$ is 
equivalent to $F$.

\thm[\cite{Matsuzaki2021}]{\label{Thm:R-moves}}{
	Two spatial surfaces are equivalent if and only if
	their diagrams are related by a finite sequence of ${\rm R}2$, ${\rm R}3$, ${\rm R}5$ and ${\rm R}6$ moves depicted in Fig.~\ref{fig Reidemeister moves for spatial surface diagrams} and isotopies in $S^2 = \mathbb{R}^{2} \cup \{ \infty \}$.   
}
\upshape

\begin{figure}[h]
  \centering
  \begin{tikzpicture}[line width=1.6pt, use Hobby shortcut]
    \begin{knot}[
      consider self intersections=true,
      ignore endpoint intersections=false,
      flip crossing/.list={},
      only when rendering/.style={}
    ]
    \strand (0,2)--(0,0);
    \strand (1,2)--(1,0);
    \strand (3,2)--(3.5,1.5)--(4,1)--(3.5,0.5)--(3,0);
    \strand (4,2)--(3.5,1.5)--(3,1)--(3.5,0.5)--(4,0);
    \end{knot}
    \draw[<->] (1.5,1)--(2.5,1);
    \draw (2,0.3) node{R$2$};
   \end{tikzpicture}
  \quad
  \begin{tikzpicture}[line width=1.6pt, use Hobby shortcut]
    \begin{knot}[
      consider self intersections=true,
        ignore endpoint intersections=false,
        flip crossing/.list={},
        only when rendering/.style={}
      ]
      \strand(2,2)--(0,0);
      \strand(1,2)--(0.5,1.5)--(0,1)--(0.5,0.5)--(1,0);
      \strand(0,2)--(2,0);
      \strand(6,2)--(4,0);
      \strand(5,2)--(6,1)--(5,0);
      \strand(4,2)--(6,0);
      \end{knot}
      \draw[<->] (2.5,1)--(3.5,1);
      \draw (3,0.3) node{R$3$};
    \end{tikzpicture}

    \vspace{5mm}
    \begin{tikzpicture}[line width=1.6pt, use Hobby shortcut]
      \begin{knot}[
        consider self intersections=true,
        ignore endpoint intersections=false,
        flip crossing/.list={},
        only when rendering/.style={}
      ]
      \strand(0.75,1.5)--(0.75,0);
      \strand(0,1.5)--(1.5,0.75)--(0,0);
      \strand(1.5,0.75)--(3,0.75);
      \strand(4.5,1.5)--(6,0.75)--(4.5,0);
      \strand(6.75,1.5)--(6.75,0);
      \strand(6,0.75)--(7.5,0.75);
      \end{knot}
      \filldraw (1.5,0.75) circle(0.15);
      \draw[<->] (3.375,0.75)--(4.125,0.75);
      \filldraw (6,0.75) circle(0.15);
      \draw (3.75,0.225) node{R$5$};
    \end{tikzpicture}
		\quad
		\begin{tikzpicture}[line width=1.6pt, use Hobby shortcut]
      \begin{knot}[
        consider self intersections=true,
        ignore endpoint intersections=false,
        flip crossing/.list={},
        only when rendering/.style={}
      ]
      \strand(0,1.5)--(1.5,0.75)--(0,0);
      \strand(0.75,1.5)--(0.75,0);
      \strand(1.5,0.75)--(3,0.75);
      \strand(4.5,1.5)--(6,0.75)--(4.5,0);
      \strand(6,0.75)--(7.5,0.75);
      \strand(6.75,1.5)--(6.75,0);
      \end{knot}
      \filldraw (1.5,0.75) circle(0.15);
      \draw[<->] (3.375,0.75)--(4.125,0.75);
      \filldraw (6,0.75) circle(0.15);
      \draw (3.75,0.225) node{R$5$};
    \end{tikzpicture}

        \vspace{5mm}
    \begin{tikzpicture}[line width=1.6pt, use Hobby shortcut]
      \begin{knot}[
        consider self intersections=true,
        ignore endpoint intersections=false,
        flip crossing/.list={},
        only when rendering/.style={}
      ]
      \strand(0,2)--(1,1.5)--(2,2);
      \strand(1,1.5)--(1,0.5);
      \strand(0,0)--(1,0.5)--(2,0);
      \strand(6,2)--(6.5,1)--(6,0);
      \strand(6.5,1)--(7.5,1);
      \strand(8,0)--(7.5,1)--(8,2);
      \end{knot}
      \filldraw (1,1.5) circle(0.15);
      \filldraw (1,0.5) circle(0.15);
      \draw[<->] (3,1)--(5,1);
      \draw (4,0.3) node{R$6$};
      \filldraw (6.5,1) circle(0.15);
      \filldraw (7.5,1) circle(0.15);
    \end{tikzpicture}
    \caption{Reidemeister moves for diagrams of spatial surfaces}{\label{fig Reidemeister moves for spatial surface diagrams}}
\end{figure}  

A \textit{Y-orientation} of a spatial trivalent graph is an assignment of orientations to \textcolor{black}{the edges of the spatial trivalent graph} 
such that no sinks or sources exist (Fig.~3).
A \textit{Y-oriented spatial trivalent graph} is 
a spatial trivalent graph with a Y-orientation.
Any spatial trivalent graph admits at least one Y-orientation.
\textcolor{black}{A Y-orientation of a diagram of a spatial trivalent graph is also defined in the same manner.}  
A \textit{Y-oriented diagram} of a spatial surface is a diagram of the spatial surface with a Y-orientation.

\begin{figure}[h]
  \centering
  \begin{tikzpicture}[line width=2pt]
    \draw (0,2)--(1,1)--(1,0);
    \draw (2,2)--(1,1);
    \draw (0.5,1.5) node{\rotatebox{315}{$\blacktriangleright$}};
    \draw (1.5,1.5) node{\rotatebox{225}{$\blacktriangleright$}};
    \draw (1,0.5) node{\rotatebox{270}{$\blacktriangleright$}};
    \filldraw (1,1) circle(0.125);
    \draw (1,-0.5) node{Y-orientation};
  \end{tikzpicture}
  \quad 
  \begin{tikzpicture}[line width=2pt]
    \draw (0,2)--(1,1)--(1,0);
    \draw (2,2)--(1,1);
    \draw (0.5,1.5) node{\rotatebox{135}{$\blacktriangleright$}};
    \draw (1.5,1.5) node{\rotatebox{45}{$\blacktriangleright$}};
    \draw (1,0.5) node{\rotatebox{90}{$\blacktriangleright$}};
    \filldraw (1,1) circle(0.125);
    \draw (1,-0.5) node{Y-orientation};
  \end{tikzpicture}
  \quad
  \begin{tikzpicture}[line width=2pt]
    \draw (0,2)--(1,1)--(1,0);
    \draw (2,2)--(1,1);
    \draw (0.5,1.5) node{\rotatebox{315}{$\blacktriangleright$}};
    \draw (1.5,1.5) node{\rotatebox{235}{$\blacktriangleright$}};
    \draw (1,0.5) node{\rotatebox{90}{$\blacktriangleright$}};
    \filldraw (1,1) circle(0.125);
    \draw (1,-0.5) node{sink};
  \end{tikzpicture}
  \quad
  \begin{tikzpicture}[line width=2pt]
    \draw (0,2)--(1,1)--(1,0);
    \draw (2,2)--(1,1);
    \draw (0.5,1.5) node{\rotatebox{135}{$\blacktriangleright$}};
    \draw (1.5,1.5) node{\rotatebox{45}{$\blacktriangleright$}};
    \draw (1,0.5) node{\rotatebox{270}{$\blacktriangleright$}};
    \filldraw (1,1) circle(0.125);
    \draw (1,-0.5) node{source};
  \end{tikzpicture}
  \caption{All orientations around trivalent vertices}{\label{fig All probable Y-orientations around trivalent vertices}}
\end{figure}

\textit{Y-oriented ${\rm R}2$, ${\rm R}3$, ${\rm R}5$, and ${\rm R}6$ moves} are \textcolor{black}{local} moves on Y-oriented diagrams whose underlying moves are 
${\rm R}2$, ${\rm R}3$, ${\rm R}5$ and ${\rm R}6$ moves depicted in Fig.~\ref{fig Reidemeister moves for spatial surface diagrams}, respectively. 
We refer to them as \textit{Y-oriented Reidemeister moves}. 

\lem{\label{Lemma:generating set}}{
Y-oriented Reidemeister moves are generated by oriented ${\rm R}2$, ${\rm R}3$ moves, 
and Y-oriented ${\rm R}5$ and ${\rm R}6$ moves depicted in Fig.~\ref{fig generating set of Y-oriented Reidemeister moves for oriented spatial surface diagrams}.
\begin{figure}[h]
  \centering
  \begin{tikzpicture}[line width=1.6pt,scale=0.875, use Hobby shortcut]
    \begin{knot}[
      consider self intersections=true,
      ignore endpoint intersections=false,
      flip crossing/.list={},
      only when rendering/.style={}
    ]
    \strand (0.4,2)--(0.4,-0.4);
    \strand (0,1.6)--(1.6,0.8)--(0,0);

    \strand (7.2,2)--(7.2,-0.4);
    \strand (6.4,0.8)--(8,0.8);
    \end{knot}
    \draw (1.6,0.8)--(3.2,0.8);
    \filldraw (1.6,0.8) circle(0.105);
    \draw (1,1.1) node{\rotatebox{-26.5}{$\blacktriangleright$}};
    \draw (1,0.5) node{\rotatebox{26.5}{$\blacktriangleright$}};
    \draw (0.4,2) node{$\blacktriangle$};
    \draw (3.2,0.8) node{$\blacktriangleright$};

    \draw[<->] (3.6,0.8)--(4.4,0.8);
    \draw (4,0.56) node[below]{R$5$(A)};

    \draw (4.8,1.6)--(6.4,0.8)--(4.8,0);
    \filldraw (6.4,0.8) circle(0.105);
    \draw (5.8,1.1) node{\rotatebox{-26.5}{$\blacktriangleright$}};
    \draw (5.8,0.5) node{\rotatebox{26.5}{$\blacktriangleright$}};
    \draw (7.2,2) node{$\blacktriangle$};
    \draw (8,0.8) node{$\blacktriangleright$};
  \end{tikzpicture}
  \begin{tikzpicture}[line width=1.6pt, scale=0.7, use Hobby shortcut]
    \begin{knot}[
      consider self intersections=true,
      ignore endpoint intersections=false,
      flip crossing/.list={},
      only when rendering/.style={}
    ]
    \strand (0,2)--(2,1)--(0,0);
    \strand (0.5,2.5)--(0.5,-0.5);

    \strand (8,1)--(10,1);
    \strand (9,2.5)--(9,-0.5);
    \end{knot}
    \draw (2,1)--(4,1);
    \filldraw (2,1) circle(0.15);
    \draw (0,2) node{\rotatebox{-26.5}{$\blacktriangleleft$}};
    \draw (0,0) node{\rotatebox{26.5}{$\blacktriangleleft$}};
    \draw (0.5,2.5) node{$\blacktriangle$};
    \draw (3,1) node{$\blacktriangleleft$};

    \draw[<->] (4.5,1)--(5.5,1);
    \draw (5,0.7) node[below]{R$5$(B)};

    \draw (6,2)--(8,1)--(6,0);
    \filldraw (8,1) circle(0.15);
    \draw (6,2) node{\rotatebox{-26.5}{$\blacktriangleleft$}};
    \draw (6,0) node{\rotatebox{26.5}{$\blacktriangleleft$}};
    \draw (9,2.5) node{$\blacktriangle$};
    \draw (9.5,1) node{$\blacktriangleleft$};
  \end{tikzpicture}
  
  \vspace*{5mm}
  \begin{tikzpicture}[line width=1.6pt,scale=0.7, use Hobby shortcut]
    \begin{knot}[
      consider self intersections=true,
      ignore endpoint intersections=false,
      flip crossing/.list={},
      only when rendering/.style={}
    ]
    \strand (0.5,2.5)--(0.5,-0.5);
    \strand (0,2)--(2,1)--(0,0);

    \strand (9,2.5)--(9,-0.5);
    \strand (8,1)--(10,1);
    \end{knot}
    \draw (2,1)--(4,1);
    \filldraw (2,1) circle(0.15);
    \draw (1.25,1.375) node{\rotatebox{-26.5}{$\blacktriangleright$}};
    \draw (1.25,0.625) node{\rotatebox{26.5}{$\blacktriangleright$}};
    \draw (0.5,-0.5) node{$\blacktriangledown$};
    \draw (4,1) node{$\blacktriangleright$};

    \draw[<->] (4.5,1)--(5.5,1);
    \draw (5,0.7) node[below]{R$5$(C)};

    \draw (6,2)--(8,1)--(6,0);
    \filldraw (8,1) circle(0.15);
    \draw (7.25,1.375) node{\rotatebox{-26.5}{$\blacktriangleright$}};
    \draw (7.25,0.625) node{\rotatebox{26.5}{$\blacktriangleright$}};
    \draw (9,-0.5) node{$\blacktriangledown$};
    \draw (10,1) node{$\blacktriangleright$};
  \end{tikzpicture}
  \begin{tikzpicture}[line width=1.6pt,scale=0.7, use Hobby shortcut]
    \begin{knot}[
      consider self intersections=true,
      ignore endpoint intersections=false,
      flip crossing/.list={},
      only when rendering/.style={}
    ]
    \strand (0,2)--(2,1)--(0,0);
    \strand (0.5,2.5)--(0.5,-0.5);

    \strand (8,1)--(10,1);
    \strand (9,2.5)--(9,-0.5);
    \end{knot}
    \draw (2,1)--(4,1);
    \filldraw (2,1) circle(0.15);
    \draw (0,2) node{\rotatebox{-26.5}{$\blacktriangleleft$}};
    \draw (0,0) node{\rotatebox{26.5}{$\blacktriangleleft$}};
    \draw (0.5,-0.5) node{$\blacktriangledown$};
    \draw (3,1) node{$\blacktriangleleft$};

    \draw[<->] (4.5,1)--(5.5,1);
    \draw (5,0.7) node[below]{R$5$(D)};

    \draw (6,2)--(8,1)--(6,0);
    \filldraw (8,1) circle(0.15);
    \draw (6,2) node{\rotatebox{-26.5}{$\blacktriangleleft$}};
    \draw (6,0) node{\rotatebox{26.5}{$\blacktriangleleft$}};
    \draw (9,-0.5) node{$\blacktriangledown$};
    \draw (9.5,1) node{$\blacktriangleleft$};
  \end{tikzpicture}
  
\vspace{2.5mm}
      \begin{tikzpicture}[line width=1.6pt]
        \draw (0,3)--(1,2)--(2,3);
        \draw (1,2)--(1,1);
        \draw (0,0)--(1,1)--(2,0);
        \draw (0.5,2.5) node{\rotatebox{-45}{$\blacktriangleright$}};
        \draw (1.5,2.5) node{\rotatebox{45}{$\blacktriangleleft$}};
        \draw (1,1.5) node{$\blacktriangledown$};
        \draw (0.5,0.5) node{\rotatebox{45}{$\blacktriangleright$}};
        \draw (1.5,0.5) node{\rotatebox{-45}{$\blacktriangleright$}};
        \filldraw (1,2) circle(0.105);
        \filldraw (1,1) circle(0.105);
    
        \draw[<->] (2.5,1.5)--(3.5,1.5);
        \draw (3,1) node{R$6$(A)};
    
        \draw (4,2.5)--(5,1.5)--(6,1.5)--(7,2.5);
        \draw (4,0.5)--(5,1.5);
        \draw (6,1.5)--(7,0.5);
        \filldraw (5,1.5) circle(0.105);
        \filldraw (6,1.5) circle(0.105);
        \draw (4.5,2) node{\rotatebox{-45}{$\blacktriangleright$}};
        \draw (4.5,1) node{\rotatebox{45}{$\blacktriangleright$}};
        \draw (5.5,1.5) node{$\blacktriangleright$};
        \draw (6.5,2) node{\rotatebox{45}{$\blacktriangleleft$}};
        \draw (6.5,1) node{\rotatebox{-45}{$\blacktriangleright$}};
      \end{tikzpicture}
      \begin{tikzpicture}[line width=1.6pt]
        \draw (0,3)--(1,2)--(2,3);
        \draw (1,2)--(1,1);
        \draw (0,0)--(1,1)--(2,0);
        \draw (0.5,2.5) node{\rotatebox{-45}{$\blacktriangleright$}};
        \draw (1.5,2.5) node{\rotatebox{45}{$\blacktriangleright$}};
        \draw (1,1.5) node{$\blacktriangledown$};
        \draw (0.5,0.5) node{\rotatebox{45}{$\blacktriangleleft$}};
        \draw (1.5,0.5) node{\rotatebox{-45}{$\blacktriangleright$}};
        \filldraw (1,2) circle(0.105);
        \filldraw (1,1) circle(0.105);
    
        \draw[<->] (2.5,1.5)--(3.5,1.5);
        \draw (3,1) node{R$6$(B)};
    
        \draw (4,2.5)--(5,1.5)--(6,1.5)--(7,2.5);
        \draw (4,0.5)--(5,1.5);
        \draw (6,1.5)--(7,0.5);
        \filldraw (5,1.5) circle(0.105);
        \filldraw (6,1.5) circle(0.105);
        \draw (4.5,2) node{\rotatebox{-45}{$\blacktriangleright$}};
        \draw (4.5,1) node{\rotatebox{45}{$\blacktriangleleft$}};
        \draw (5.5,1.5) node{$\blacktriangleright$};
        \draw (6.5,2) node{\rotatebox{45}{$\blacktriangleright$}};
        \draw (6.5,1) node{\rotatebox{-45}{$\blacktriangleright$}};
      \end{tikzpicture}
      \vspace{2.5mm}
      \begin{tikzpicture}[line width=1.6pt]
        \draw (-2,2.5)--(-3,1.5)--(-4,1.5)--(-5,2.5);
        \draw (-5,0.5)--(-4,1.5);
        \draw (-2,0.5)--(-3,1.5);
        \filldraw (-3,1.5) circle(0.105);
        \filldraw (-4,1.5) circle(0.105);
        \draw (-4.5,2) node{\rotatebox{-45}{$\blacktriangleright$}};
        \draw (-4.5,1) node{\rotatebox{45}{$\blacktriangleleft$}};
        \draw (-3.5,1.5) node{$\blacktriangleleft$};
        \draw (-2.5,2) node{\rotatebox{45}{$\blacktriangleleft$}};
        \draw (-2.5,1) node{\rotatebox{-45}{$\blacktriangleright$}};
    
        \draw[<->] (-0.5,1.5)--(-1.5,1.5);
        \draw (-1,1) node{R$6$(C)}; 
    
        \draw (0,3)--(1,2)--(2,3);
        \draw (1,2)--(1,1);
        \draw (0,0)--(1,1)--(2,0);
        \draw (0.5,2.5) node{\rotatebox{-45}{$\blacktriangleright$}};
        \draw (1.5,2.5) node{\rotatebox{45}{$\blacktriangleleft$}};
        \draw (1,1.5) node{$\blacktriangledown$};
        \draw (0.5,0.5) node{\rotatebox{45}{$\blacktriangleleft$}};
        \draw (1.5,0.5) node{\rotatebox{-45}{$\blacktriangleright$}};
        \filldraw (1,2) circle(0.105);
        \filldraw (1,1) circle(0.105);
    
        \draw[<->] (2.5,1.5)--(3.5,1.5);
        \draw (3,1) node{R$6$(D)};
    
        \draw (4,2.5)--(5,1.5)--(6,1.5)--(7,2.5);
        \draw (4,0.5)--(5,1.5);
        \draw (6,1.5)--(7,0.5);
        \filldraw (5,1.5) circle(0.105);
        \filldraw (6,1.5) circle(0.105);
        \draw (4.5,2) node{\rotatebox{-45}{$\blacktriangleright$}};
        \draw (4.5,1) node{\rotatebox{45}{$\blacktriangleleft$}};
        \draw (5.5,1.5) node{$\blacktriangleright$};
        \draw (6.5,2) node{\rotatebox{45}{$\blacktriangleleft$}};
        \draw (6.5,1) node{\rotatebox{-45}{$\blacktriangleright$}};
      \end{tikzpicture}
      \vspace{2.5mm}
      \begin{tikzpicture}[line width=1.6pt]
        \draw (-2,2.5)--(-3,1.5)--(-4,1.5)--(-5,2.5);
        \draw (-5,0.5)--(-4,1.5);
        \draw (-2,0.5)--(-3,1.5);
        \filldraw (-3,1.5) circle(0.105);
        \filldraw (-4,1.5) circle(0.105);
        \draw (-4.5,2) node{\rotatebox{-45}{$\blacktriangleright$}};
        \draw (-4.5,1) node{\rotatebox{45}{$\blacktriangleleft$}};
        \draw (-3.5,1.5) node{$\blacktriangleleft$};
        \draw (-2.5,2) node{\rotatebox{45}{$\blacktriangleright$}};
        \draw (-2.5,1) node{\rotatebox{-45}{$\blacktriangleleft$}};
    
        \draw[<->] (-0.5,1.5)--(-1.5,1.5);
        \draw (-1,1) node{R$6$(E)}; 
    
        \draw (0,3)--(1,2)--(2,3);
        \draw (1,2)--(1,1);
        \draw (0,0)--(1,1)--(2,0);
        \draw (0.5,2.5) node{\rotatebox{-45}{$\blacktriangleright$}};
        \draw (1.5,2.5) node{\rotatebox{45}{$\blacktriangleright$}};
        \draw (1,1.5) node{$\blacktriangledown$};
        \draw (0.5,0.5) node{\rotatebox{45}{$\blacktriangleleft$}};
        \draw (1.5,0.5) node{\rotatebox{-45}{$\blacktriangleleft$}};
        \filldraw (1,2) circle(0.105);
        \filldraw (1,1) circle(0.105);
    
        \draw[<->] (2.5,1.5)--(3.5,1.5);
        \draw (3,1) node{R$6$(F)};
    
        \draw (4,2.5)--(5,1.5)--(6,1.5)--(7,2.5);
        \draw (4,0.5)--(5,1.5);
        \draw (6,1.5)--(7,0.5);
        \filldraw (5,1.5) circle(0.105);
        \filldraw (6,1.5) circle(0.105);
        \draw (4.5,2) node{\rotatebox{-45}{$\blacktriangleright$}};
        \draw (4.5,1) node{\rotatebox{45}{$\blacktriangleleft$}};
        \draw (5.5,1.5) node{$\blacktriangleright$};
        \draw (6.5,2) node{\rotatebox{45}{$\blacktriangleright$}};
        \draw (6.5,1) node{\rotatebox{-45}{$\blacktriangleleft$}};
      \end{tikzpicture}
  \caption{Y-oriented R$5$ and R$6$ moves}{\label{fig generating set of Y-oriented Reidemeister moves for oriented spatial surface diagrams}}
\end{figure}}\upshape

\begin{proof}
  A Y-oriented R$5$ move as illustrated in Fig.~\ref{Fig:Y-oriented_R5-move} is realized through a sequence of Y-oriented Reidemeister moves as depicted in Fig.~\ref{Fig:Seq}.
  \begin{figure}[h]
    \centering
    \begin{tikzpicture}[line width=1.6pt,scale=0.875, use Hobby shortcut]
      \begin{knot}[
        consider self intersections=true,
        ignore endpoint intersections=false,
        flip crossing/.list={},
        only when rendering/.style={}
      ]
      \strand (0.4,2)--(0.4,-0.4);
      \strand (0,1.6)--(1.6,0.8)--(0,0);
  
      \strand (7.2,2)--(7.2,-0.4);
      \strand (6.4,0.8)--(8,0.8);
      \end{knot}
      \draw (1.6,0.8)--(3.2,0.8);
      \filldraw (1.6,0.8) circle(0.105);
      \draw (1,1.1) node{\rotatebox{-26.5}{$\blacktriangleright$}};
      \draw (1,0.5) node{\rotatebox{26.5}{$\blacktriangleleft$}};
      \draw (0.4,2) node{$\blacktriangle$};
      \draw (3.2,0.8) node{$\blacktriangleright$};
  
      \draw[<->] (3.6,0.8)--(4.4,0.8);
      \draw (4,0.56) node[below]{R$5$};
  
      \draw (4.8,1.6)--(6.4,0.8)--(4.8,0);
      \filldraw (6.4,0.8) circle(0.105);
      \draw (5.8,1.1) node{\rotatebox{-26.5}{$\blacktriangleright$}};
      \draw (5.8,0.5) node{\rotatebox{26.5}{$\blacktriangleleft$}};
      \draw (7.2,2) node{$\blacktriangle$};
      \draw (8,0.8) node{$\blacktriangleright$};
    \end{tikzpicture}
    \caption{A Y-oriented R$5$ move}{\label{Fig:Y-oriented_R5-move}}
  \end{figure}
  Similarly,
  other Y-oriented R$5$ moves not included in Fig.~\ref{fig generating set of Y-oriented Reidemeister moves for oriented spatial surface diagrams}
  can be realized by a sequence of moves consisting of R$5$(A) -- R$5$(D) and oriented R$2$ moves.

  \begin{figure}[h]
    \centering
    \begin{tikzpicture}[line width=1.6pt,scale=0.875, use Hobby shortcut]
      \begin{knot}[
        consider self intersections=true,
        ignore endpoint intersections=false,
        flip crossing/.list={},
        only when rendering/.style={}
      ]
      \strand (0.4,2)--(0.4,-0.4);
      \strand (0,1.6)--(1.6,0.8)--(0,0);
  
      \strand (5.2,2) .. (5.2,1.8) .. (6.8,0.2) .. (5.2,0.7) .. (5.2,-0.2) .. (5.2,-0.4);
      \strand (4.8,1.6)--(6.4,0.8)--(4.8,0);
      \strand (6.4,0.8)--(8,0.8);

      \strand (12,2)--(12,-0.4);
      \strand (11.2,0.8) -- (12.8,0.8);
      \end{knot}
      \draw (1.6,0.8)--(3.2,0.8);
      \filldraw (1.6,0.8) circle(0.105);
      \draw (1,1.1) node{\rotatebox{-26.5}{$\blacktriangleright$}};
      \draw (1,0.5) node{\rotatebox{26.5}{$\blacktriangleleft$}};
      \draw (0.4,2) node{$\blacktriangle$};
      \draw (3.2,0.8) node{$\blacktriangleright$};
  
      \draw[<->] (3.6,0.8)--(4.4,0.8);
      \draw (4,0.56) node[below]{R$5$(A)};
  
      \filldraw (6.4,0.8) circle(0.105);
      \draw (5.8,1.1) node{\rotatebox{-26.5}{$\blacktriangleright$}};
      \draw (5.4,0.3) node{\rotatebox{26.5}{$\blacktriangleleft$}};
      \draw (5.2,2) node{$\blacktriangle$};
      \draw (8,0.8) node{$\blacktriangleright$};
      \draw (6.4,0.8) -- (7,0.8);

      \draw[<->] (8.4,0.8)--(9.2,0.8);
      \draw (8.8,0.56) node[below]{R$2$};
  
      \filldraw (11.2,0.8) circle(0.105);
      \draw (10.6,1.1) node{\rotatebox{-26.5}{$\blacktriangleright$}};
      \draw (10.6,0.5) node{\rotatebox{26.5}{$\blacktriangleleft$}};
      \draw (12,2) node{$\blacktriangle$};
      \draw (12.8,0.8) node{$\blacktriangleright$};
      \draw (9.6,1.6)--(11.2,0.8)--(9.6,0);
    \end{tikzpicture}
    \caption{A sequence of Y-oriented Reidemeister moves}{\label{Fig:Seq}}
  \end{figure}

  In Fig.~\ref{fig generating set of Y-oriented Reidemeister moves for oriented spatial surface diagrams},
  all Y-oriented R$6$ moves are listed.
  Thus Y-oriented Reidemeister moves are generated by oriented R$2$, R$3$ moves and Y-oriented R$5$, R$6$ moves depicted in Fig.~\ref{fig generating set of Y-oriented Reidemeister moves for oriented spatial surface diagrams}.
\end{proof}

Let $D$ and $D'$ be Y-oriented diagrams such that
$D'$ is obtained from \textcolor{black}{$D$ by} a Y-oriented Reidemeister move. 
Let $U$ be a disk in $\mathbb{R}^{2}$ such that $D \setminus U = D' \setminus U$ and
$D \cap U$ and $D' \cap U$ are as in Fig.~\ref{fig Reidemeister moves for spatial surface diagrams} 
when we forget the Y-orientations.
We call $U$ \textit{a support of the Y-oriented Reidemeister move} and
we say that $D'$ \textit{is obtained from $D$ by a Y-oriented Reidemeister move with support $U$}.



In this paper, we call an operation that reverses an orientation of an $S^{1}$-component an \textit{inverse move}.




\textcolor{black}{Two Y-oriented diagrams $D$ and $D'$, 
whose underlying diagrams are the same, 
are related by a finite sequence of Y-oriented Reidemeister moves, inverse moves, and isotopies in $S^{2}$ (see \cite{Ishii2015}, \cite{Matsuzaki-Murao2023}).
In combination with Theorem~\ref{Thm:R-moves}, 
a Reidemeister-type theorem on Y-oriented diagrams is introduced in \cite{Matsuzaki-Murao2023}.}

\thm[\cite{Matsuzaki-Murao2023}]{\label{thm Y-oriented R-moves}}{
  Two spatial surfaces are equivalent if and only if 
  their Y-oriented diagrams are related by a finite sequence of Y-oriented Reidemeister moves, 
  isotopies in $S^2$, 
   and inverse moves.
}

\upshape


\rem{A spatial surface is said to be \textit{annular} if every connected component of the spatial surface is homeomorphic to an annulus.
An annular spatial surface corresponds to a framed link in $S^{3}$.
In this sense,
a spatial surface is a generalization of a framed link in $S^{3}$.

Diagrams of annular spatial surfaces are link diagrams.
Conversely,
link diagrams \textcolor{black}{correspond} to annular spatial surfaces.
According to Theorem~\ref{Thm:R-moves},
any two diagrams presenting equivalent annular spatial surfaces are related by a finite sequence of R$2$ moves, R$3$ moves, and isotopies of $S^{2}$.

Y-oriented diagrams of annular spatial surfaces are oriented link diagrams.
From Theorem~\ref{thm Y-oriented R-moves}, 
any two oriented link diagrams presenting the same annular spatial surface 
are related by a finite sequence of oriented R$2$ moves, oriented R$3$ moves, isotopies on $S^{2}$, and inverse moves.
}

\section{Colorings for diagrams of spatial surfaces}

Let $D$ be a Y-oriented diagram of a spatial surface $F$.
An \textit{arc} of $D$ means a simple arc or a simple loop which is obtained from $D$ by cutting the diagram at under crossings and vertices.
We denote the set of all arcs by $\mathcal{A}(D)$.
\defi{\label{Def:Coloring}}{
Let $X$ be a groupoid rack.
An \textit{$X$-coloring} or a \textit{coloring by $X$} of a Y-oriented diagram $D$ is a map $C : \mathcal{A}(D) \to X$ satisfying the following conditions:
\begin{itemize}
  \item[(i)] For each crossing of $D$, $C$ satisfies $C(a_{i}) \ast C(a_{j}) = C(a_{k})$,
  where $a_{i}, a_{j}, a_{k} \in \mathcal{A}(D)$ are as shown on the left side of Fig.~\ref{fig coloring conditions around crossings and trivalent vertices}.
  \item[(ii)] For each vertex of $D$, $C$ satisfies $\mathrm{cod}(C(a_{i})) = \mathrm{dom}(C(a_{j}))$ and
  $C(a_{i}) C(a_{j}) = C(a_{k})$\textcolor{black}{,} where $a_{i}, a_{j}, a_{k} \in \mathcal{A}(D)$ are as shown in the center or \textcolor{black}{the} right side of Fig.~\ref{fig coloring conditions around crossings and trivalent vertices}.
\end{itemize}
}
%
%
We denote the set of all $X$-colorings of $D$ by $\textrm{Col}_{X}(D)$.

\begin{figure}[h]
  \centering
    \begin{tikzpicture}[line width=2pt, use Hobby shortcut]
      \begin{knot}[
        consider self intersections=true,
        ignore endpoint intersections=false,
        flip crossing/.list={}
        ]
      \strand (2,2)--(0,0);
      \strand (0,2)--(2,0);
      \end{knot}  
      \draw (0,0) node{\rotatebox{225}{$\blacktriangleright$}};
      \draw (-0.3,2) node{$a_{i}$};
      \draw (2.3,2) node{$a_{j}$};
      \draw (2.3,0) node{$a_{k}$};
      \draw (1,-0.5) node{$C(a_{i}) \ast C(a_{j}) = C(a_{k})$};
    \end{tikzpicture}
    \quad
    \begin{tikzpicture}[line width=2pt]
      \draw (0,2)--(1,1)--(1,0);
      \draw (2,2)--(1,1);
      \draw (0.5,1.5) node{\rotatebox{315}{$\blacktriangleright$}};
      \draw (1.5,1.5) node{\rotatebox{225}{$\blacktriangleright$}};
      \draw (1,0.5) node{\rotatebox{270}{$\blacktriangleright$}};
      \draw (-0.3,2) node{$a_{i}$};
      \draw (2.3,2) node{$a_{j}$};
      \draw (1.5,0) node{$a_{k}$};
      \filldraw (1,1) circle(0.125);
      \draw (1,-0.5) node{$C(a_{i})C(a_{j}) = C(a_{k})$};
    \end{tikzpicture}
    \quad 
    \begin{tikzpicture}[line width=2pt]
      \draw (1,2)--(1,1)--(0,0);
      \draw (2,0)--(1,1);
      \draw (0.5,0.5) node{\rotatebox{225}{$\blacktriangleright$}};
      \draw (1.5,0.5) node{\rotatebox{315}{$\blacktriangleright$}};
      \draw (1,1.5) node{\rotatebox{270}{$\blacktriangleright$}};
      \filldraw (1,1) circle(0.125);
      \draw (1.5,2) node{$a_{k}$};
      \draw (-0.3,0) node{$a_{i}$};
      \draw (2.3,0) node{$a_{i}$};
      \draw (1,-0.5) node{$C(a_{i})C(a_{j}) = C(a_{k})$};
    \end{tikzpicture} 
    \caption{$X$-coloring conditions}{\label{fig coloring conditions around crossings and trivalent vertices}}  
\end{figure}

\thm{\label{thm coloring invariant}}{
	Let $X$ be a groupoid rack and let $D$ and $D'$ be Y-oriented diagrams of spatial surfaces $F$ and $F'$, respectively.  If $F$ and $F'$ are equivalent, then  
	there is a bijection between $\textrm{Col}_{X}(D)$ and $\textrm{Col}_{X}(D')$.
  In particular, 
  the cardinality of $\textrm{Col}_{X}(D)$ is an invariant of a spatial surface.
}
\upshape

\textcolor{black}{We will prove Theorem~\ref{thm coloring invariant} in Section~\ref{Sec:proof}.}
\rem{
  As seen in Examples \ref{ex MGR} and \ref{ex heap rack},
  multiple group racks and heap racks are regarded as groupoid racks.
  Then
  colorings using multiple group racks defined in \cite{Ishii-Matsuzaki-Murao2020} and
   colorings using heap racks defined in \cite{Saito-Zappala2024}
  are regarded as colorings using groupoid racks in our sense (Definition~\ref{Def:Coloring}).
}

\section{A proof of Theorem \ref{thm coloring invariant}}{\label{Sec:proof}}

In this section,
let $X$ be a groupoid rack with the good involution $\rho : X \to X$ sending $x$ to $x^{-1}$ as in \textcolor{black}{Proposition~\ref{Prop:Good_involution}} 
and let $P$ denote the set $\left\{(f, g) \in X^{2} \mid \mathrm{cod}(f) = \mathrm{dom}(g)\right\}$.



\lem{\label{Lemma:well-def}}{
  For any $a, b, c \in X$, the following \textcolor{black}{four conditions} are equivalent:
  \begin{enumerate}
    \item[\rm{(i)}] $(a, b) \in P$ and $ab = c$.
    \item[\rm{(ii)}] $(\rho(a), c) \in P$ and $\rho(a)c = b$.
    \item[\rm{(iii)}] $(c, \rho(b)) \in P$ and $c \rho(b) = a$.
    \item[\rm{(iv)}] $(\rho(b), \rho(a)) \in P$ and $\rho(b) \rho(a) = \rho(c)$.
  \end{enumerate}
}\upshape

\begin{proof}
  In any groupoid, 
  $$ab = c \Leftrightarrow a^{-1} c = b \Leftrightarrow c b^{-1} = a \Leftrightarrow b^{-1} a^{-1} = c^{-1}.$$
  Since $\rho(x) = x^{-1}$ for $x \in X$, 
  we see that 
  (i), (ii), (iii) and (iv) are equivalent.
\end{proof}


\lem{\label{Lemma:Coloring condition}}{
  \textcolor{black}{Let $D$ be a Y-oriented diagram $D$, $C: \mathcal{A}(D) \to X$ a map, and 
  $v$ a vertex in $D$. 
  Assume that the three arcs around $v$ are labeled $a$, $b$, and $c$ by $C$ 
  as in Fig.~\ref{Fig:Coloring_around_a_vertex}.}
  \begin{figure}[h]
    \centering
      \begin{tikzpicture}[line width=2pt]
        \draw (0,2)--(1,1)--(1,0);
        \draw (2,2)--(1,1);
        \draw (0.5,1.5) node{\rotatebox{315}{$\blacktriangleright$}};
        \draw (1.5,1.5) node{\rotatebox{225}{$\blacktriangleright$}};
        \draw (1,0.5) node{\rotatebox{270}{$\blacktriangleright$}};
        \draw (-0.3,2) node{$a$};
        \draw (2.3,2) node{$b$};
        \draw (1.5,0) node{$c$};
        \filldraw (1,1) circle(0.125);
      \end{tikzpicture}
      \quad
      \begin{tikzpicture}[line width=2pt]
        \draw (1,2) -- (1,1);
        \draw (0,0) -- (1,1) -- (2,0);
        \draw (1,1.5) node{$\blacktriangledown$};
        \draw (0.5,0.5) node{\rotatebox{45}{$\blacktriangleleft$}};
        \draw (1.5,0.5) node{\rotatebox{-45}{$\blacktriangleright$}};
        \draw (-0.3,0) node{$a$};
        \draw (2.3,0) node{$b$};
        \draw (1.3,2) node{$c$};
         \filldraw (1,1) circle(0.125);
      \end{tikzpicture}
      \caption{An $X$-coloring $C$ around $v$}{\label{Fig:Coloring_around_a_vertex}}
    \end{figure}
    \textcolor{black}{Then $C$ satisfies the condition of an $X$-coloring at $v$ if and only if $a$, $b$, and $c$ satisfy the conditions $\mathrm{(i)}$, $\mathrm{(ii)}$, $\mathrm{(iii)}$ and $\mathrm{(iv)}$ of Lemma~\ref{Lemma:well-def}.}
}\upshape

\begin{proof}
  The condition of an $X$-coloring at vertex $v$ is the condition (i).
  By Lemma~\ref{Lemma:well-def}, 
  we \textcolor{black}{obtain the claim}.
\end{proof}

\lem{\label{Lemma:aaa}}{
  For any $a, b, c, x \in X$,
  the following three statements are equivalent.
  \begin{itemize}
    \item[\rm{(i)}] $(a, b) \in P$ and $ab = c$.
    \item[\rm{(ii)}] $(a \ast x, b \ast x) \in P$ and $(a \ast x) (b \ast x) = c \ast x$.
    \item[\rm{(iii)}] $(a \ast \rho(x), b \ast \rho(x)) \in P$ and $(a \ast \rho(x)) (b \ast \rho(x)) = c \ast \rho(x)$.
  \end{itemize}
}\upshape

\begin{proof}

  $\mathrm{(i)} \Rightarrow \mathrm{(ii)}$: 
  Suppose \textcolor{black}{that} $(a, b) \in P$ and $ab = c$.
  By the condition (iii) of Definition~\ref{def:groupoindrack},
  we have 
  $$(a \ast x, b \ast x) \in P \ \mbox{and}\  (a \ast x) (b \ast x) = (ab) \ast x = c \ast x.$$ 

  $\mathrm{(ii)} \Rightarrow \mathrm{(i)}$: 
  Suppose \textcolor{black}{that} $(a \ast x, b \ast x) \in P$ and $(a \ast x)(b \ast x) = c \ast x$.
  Since $\rho$ is a good involution of the rack $X$, 
  by the condition (iii) of Definition~\ref{def:groupoindrack},
  we have 
  $$((a \ast x) \ast \rho(x), (b \ast x) \ast \rho(x)) = (a, b) \in P$$ and 
  $$ab = ((a \ast x) \ast \rho(x))((b \ast x) \ast \rho(x)) 
  = ((a \ast x)(b \ast x)) \ast \rho(x) = (c \ast x) \ast \rho(x) = c.$$

  The equivalence between (i) and (iii) is similarly ensured by the condition (iii) of Definition~\ref{def:groupoindrack}. 
\end{proof}

\lem{\label{Lemma:Invariance_under_R5-moves}}{
  Let $D$ and $D'$ be Y-oriented diagrams of spatial surfaces.
  If $D'$ is obtained from $D$ by a Y-oriented R$5$ move with support $U$,
  there is a bijection between $\mathrm{Col}_{X}(D)$ and $\mathrm{Col}_{X}(D')$.
}\upshape

\begin{proof}
Suppose that $D \cap U$ (resp. $D' \cap U$) corresponds to the left (resp. right) side of the Y-oriented R$5$(A) move depicted in Fig.~\ref{fig generating set of Y-oriented Reidemeister moves for oriented spatial surface diagrams}.

We see that for any $X$-coloring $C$ of $D$, the restriction of $C$ to $D \cap U$ \textcolor{black}{can be assumed to be} 
as illustrated on the left side of the Y-oriented R$5$(A) move depicted in Fig.~\ref{fig Y-oriented R$5$ moves},
where $a, b, x \in X$ such that $(a \ast^{-1} x, b \ast^{-1} x) \in P$.
\textcolor{black}{Since $\rho$ is a good involution of $X$,}
by Lemma~\ref{Lemma:aaa},
we have $$(a, b) \in P\ \mbox{and}\ (a \ast^{-1} x)(b \ast^{-1} x) = (a \ast \rho(x)) (b \ast \rho(x)) = (ab) \ast \rho(x) = (ab) \ast^{-1} x.$$
Since $(a, b) \in P$ and $(a \ast^{-1} x)(b \ast^{-1} x) = (ab) \ast^{-1} x$,
we see that the restriction of $C$ to $D \setminus U (= D' \setminus U)$ uniquely extends to an $X$-coloring $C'$ of $D'$.

Conversely,
for any $X$-coloring $C'$ of $D'$,
the restriction of $C'$ to $D' \cap U$ can be assumed to be 
as illustrated on the right side of the Y-oriented R$5$(A) move depicted in Fig.~\ref{fig Y-oriented R$5$ moves},
where $a, b, x \in X$ such that $(a, b) \in P$.
\textcolor{black}{Since $\rho$ is a good involution of $X$, }
by Lemma~\ref{Lemma:aaa},
we have 
$$(a \ast^{-1} x, b \ast^{-1} x) \in P\ \mbox{and}\ (a \ast^{-1} x)(b \ast^{-1} x) = (a \ast \rho(x))(b \ast \rho(x)) = (ab) \ast \rho(x) = (ab) \ast^{-1} x.$$
Thus,
we see that the restriction of $C'$ to $D' \setminus U (= D \setminus U)$ uniquely extends to an $X$-coloring $C$ of $D$.

Therefore,
the map $\mathrm{Col}_{X}(D) \to \mathrm{Col}_{X}(D')$ which sends $C$ to $C'$ as above is a bijection.

When $D$ and $D'$ are related by one of \textcolor{black}{the} other Y-oriented R$5$ moves depicted in Fig.~\ref{fig generating set of Y-oriented Reidemeister moves for oriented spatial surface diagrams},
the existence of a bijection between $\mathrm{Col}_{X}(D)$ and $\mathrm{Col}_{X}(D')$ is similarly ensured 
by Lemmas~\ref{Lemma:Coloring condition}, \ref{Lemma:aaa}, and 
$$(a \ast x) (b \ast x) = (ab) \ast x \quad (a,b, x \in X)$$
in the case of a Y-oriented R$5$(C) move, and
$$(x \ast a) \ast b = x \ast (ab) \quad (a, b, x \in X)$$
in the case of a Y-oriented R$5$(B) move or a Y-oriented R$5$(D) move (\textcolor{black}{see} Fig.~\ref{fig Y-oriented R$5$ moves}).
\end{proof}

\begin{figure}[h]
  \centering
  \begin{tikzpicture}[line width=1.6pt, use Hobby shortcut]
    \begin{knot}[
      consider self intersections=true,
      ignore endpoint intersections=false,
      flip crossing/.list={},
      only when rendering/.style={}
    ]
    \strand (0.5,2.5)--(0.5,-0.5);
    \strand (0,2)--(2,1)--(0,0);

    \strand (9.5,2.5)--(9.5,-0.5);
    \strand (8.5,1)--(10.5,1);
    \end{knot}
    \draw (2,1)--(4,1);
    \filldraw (2,1) circle(0.15);
    \draw (1.25,1.375) node{\rotatebox{-26.5}{$\blacktriangleright$}};
    \draw (1.25,0.625) node{\rotatebox{26.5}{$\blacktriangleright$}};
    \draw (0.5,2.5) node{$\blacktriangle$};
    \draw (4,1) node{$\blacktriangleright$};
    \draw (0,2) node[left]{$b$};
    \draw (0,0) node[left]{$a$};
    \draw (0.5,2.55) node[above]{$x$};
    \draw (1.3,1.55) node[above]{$b \ast^{-1} x$};
    \draw (1.3,0.45) node[below]{$a \ast^{-1} x$};
    \draw (3.2,0.6) node{$(a \ast^{-1} x) (b \ast^{-1} x)$};

    \draw[<->] (5,1)--(6,1);
    \draw (5.5,0.5) node{R$5$(A)};

    \draw (6.5,2)--(8.5,1)--(6.5,0);
    \filldraw (8.5,1) circle(0.15);
    \draw (7.75,1.375) node{\rotatebox{-26.5}{$\blacktriangleright$}};
    \draw (7.75,0.625) node{\rotatebox{26.5}{$\blacktriangleright$}};
    \draw (9.5,2.5) node{$\blacktriangle$};
    \draw (10.5,1) node{$\blacktriangleright$};
    \draw (6.5,2) node[left]{$b$};
    \draw (6.5,0) node[left]{$a$};
    \draw (9,1.3) node{$ab$};
    \draw (9.5,2.5) node[above]{$x$};
    \draw (9.5,0.5) node[right]{$(ab) \ast^{-1} x$};
  \end{tikzpicture}
  \begin{tikzpicture}[line width=1.6pt, use Hobby shortcut]
    \begin{knot}[
      consider self intersections=true,
      ignore endpoint intersections=false,
      flip crossing/.list={},
      only when rendering/.style={}
    ]
    \strand (0,2)--(2,1)--(0,0);
    \strand (0.5,2.5)--(0.5,-0.5);

    \strand (8,1)--(10,1);
    \strand (9,2.5)--(9,-0.5);
    \end{knot}
    \draw (2,1)--(4,1);
    \filldraw (2,1) circle(0.15);
    \draw (0,2) node{\rotatebox{-26.5}{$\blacktriangleleft$}};
    \draw (0,0) node{\rotatebox{26.5}{$\blacktriangleleft$}};
    \draw (0.5,2.5) node{$\blacktriangle$};
    \draw (3,1) node{$\blacktriangleleft$};
    \draw (0,2) node[left]{$a$};
    \draw (0,0) node[left]{$b$};
    \draw (0.5,-0.5) node[below]{$(x \ast a) \ast b$};
    \draw (3,0.7) node{$ab$};
    \draw (0.5,1) node[left]{$x \ast a$};
    \draw (0.5,2.5) node[above]{$x$};

    \draw[<->] (4.5,1)--(5.5,1);
    \draw (5,0.5) node{R$5$(B)};

    \draw (6,2)--(8,1)--(6,0);
    \filldraw (8,1) circle(0.15);
    \draw (6,2) node{\rotatebox{-26.5}{$\blacktriangleleft$}};
    \draw (6,0) node{\rotatebox{26.5}{$\blacktriangleleft$}};
    \draw (9,2.5) node{$\blacktriangle$};
    \draw (9.5,1) node{$\blacktriangleleft$};
    \draw (6,2) node[left]{$a$};
    \draw (6,0) node[left]{$b$};
    \draw (9,-0.5) node[below]{$x \ast (ab)$};
    \draw (10,1) node[right]{$ab$};
    \draw (9,2.5) node[above]{$x$};
  \end{tikzpicture}

  \begin{tikzpicture}[line width=1.6pt, use Hobby shortcut]
    \begin{knot}[
      consider self intersections=true,
      ignore endpoint intersections=false,
      flip crossing/.list={},
      only when rendering/.style={}
    ]
    \strand (0.5,2.5)--(0.5,-0.5);
    \strand (0,2)--(2,1)--(0,0);

    \strand (9,2.5)--(9,-0.5);
    \strand (8,1)--(10,1);
    \end{knot}
    \draw (2,1)--(4,1);
    \filldraw (2,1) circle(0.15);
    \draw (1.25,1.375) node{\rotatebox{-26.5}{$\blacktriangleright$}};
    \draw (1.25,0.625) node{\rotatebox{26.5}{$\blacktriangleright$}};
    \draw (0.5,-0.5) node{$\blacktriangledown$};
    \draw (4,1) node{$\blacktriangleright$};
    \draw (0,2) node[left]{$b$};
    \draw (0,0) node[left]{$a$};
    \draw (0.5,2.5) node[above]{$x$};
    \draw (1.25,1.5) node[above]{$b \ast x$};
    \draw (1.25,0.5) node[below]{$a \ast x$};
    \draw (3.15,0.6) node{$(a \ast x) (b \ast x)$};

    \draw[<->] (4.5,1)--(5.5,1);
    \draw (5,0.5) node{R$5$(C)};

    \draw (6,2)--(8,1)--(6,0);
    \filldraw (8,1) circle(0.15);
    \draw (7.25,1.375) node{\rotatebox{-26.5}{$\blacktriangleright$}};
    \draw (7.25,0.625) node{\rotatebox{26.5}{$\blacktriangleright$}};
    \draw (9,-0.5) node{$\blacktriangledown$};
    \draw (10,1) node{$\blacktriangleright$};
    \draw (6,2) node[left]{$b$};
    \draw (6,0) node[left]{$a$};
    \draw (8.5,1.3) node{$ab$};
    \draw (9,2.5) node[above]{$x$};
    \draw (9.5,0.5) node[right]{$(ab) \ast x$};
  \end{tikzpicture}
  \begin{tikzpicture}[line width=1.6pt, use Hobby shortcut]
    \begin{knot}[
      consider self intersections=true,
      ignore endpoint intersections=false,
      flip crossing/.list={},
      only when rendering/.style={}
    ]
    \strand (0,2)--(2,1)--(0,0);
    \strand (0.5,2.5)--(0.5,-0.5);

    \strand (8,1)--(10,1);
    \strand (9,2.5)--(9,-0.5);
    \end{knot}
    \draw (2,1)--(4,1);
    \filldraw (2,1) circle(0.15);
    \draw (0,2) node{\rotatebox{-26.5}{$\blacktriangleleft$}};
    \draw (0,0) node{\rotatebox{26.5}{$\blacktriangleleft$}};
    \draw (0.5,-0.5) node{$\blacktriangledown$};
    \draw (3,1) node{$\blacktriangleleft$};
    \draw (0,2) node[left]{$a$};
    \draw (0,0) node[left]{$b$};
    \draw (0.5,-0.5) node[below]{$(x \ast a) \ast b$};
    \draw (3,0.7) node{$ab$};
    \draw (0.5,1) node[left]{$x \ast a$};
    \draw (0.5,2.5) node[above]{$x$};

    \draw[<->] (4.5,1)--(5.5,1);
    \draw (5,0.5) node{R$5$(D)};

    \draw (6,2)--(8,1)--(6,0);
    \filldraw (8,1) circle(0.15);
    \draw (6,2) node{\rotatebox{-26.5}{$\blacktriangleleft$}};
    \draw (6,0) node{\rotatebox{26.5}{$\blacktriangleleft$}};
    \draw (9,-0.5) node{$\blacktriangledown$};
    \draw (9.5,1) node{$\blacktriangleleft$};
    \draw (6,2) node[left]{$a$};
    \draw (6,0) node[left]{$b$};
    \draw (9,-0.5) node[below]{$x \ast (ab)$};
    \draw (9.5,0.5) node[right]{$ab$};
    \draw (9,2.5) node[above]{$x$};
  \end{tikzpicture}
  \caption{Y-oriented R$5$ moves with $X$-colorings}{\label{fig Y-oriented R$5$ moves}}
\end{figure}

\lem{\label{Lemma:Invariance_under_R6-moves}}{
  Let $D$ and $D'$ be Y-oriented diagrams of spatial surfaces.
  If $D'$ is obtained from $D$ by a Y-oriented R$6$ move with support $U$,
  there is a bijection between $\mathrm{Col}_{X}(D)$ and $\mathrm{Col}_{X}(D')$.
}

\begin{proof}
Suppose that $D \cap U$ (resp. $D' \cap U$) corresponds to the left (resp. right) side of the Y-oriented R$6$(C) move depicted in Fig.~\ref{fig generating set of Y-oriented Reidemeister moves for oriented spatial surface diagrams}.

By Lemma~\ref{Lemma:Coloring condition}, we see that for any $X$-coloring $C$ of $D$,
the restriction of $C$ to $D \cap U$ \textcolor{black}{can be assumed to be} 
as illustrated on the left side of the Y-oriented R$6$(C) move depicted in Fig.~\ref{fig Y-ori'd IH-moves ver 2},
where $a, b, c \in X$ such that
$(\rho(b), a) \in P$ and $(\rho(\rho(b)a), c) \in P$. 
By Lemma~\ref{Lemma:well-def},
\textcolor{black}{remarking that $\rho$ is involutive,}
$(\rho(b), a) \in P$ implies \textcolor{black}{that} $$(\rho(a), b) \in P\ \mbox{and}\ \rho(a)b = \rho(\rho(b)a).$$
Since $\rho(\rho(b)a) = \rho(a)b$ and $(\rho(\rho(b)a), c) \in P$, 
we have $$(\rho(a)b, c) \in P.$$
By the associativity of \textcolor{black}{the} groupoid $X$,
it follows that $$(b, c) \in P,\ (\rho(a), bc) \in P,\ \mbox{and}\ (\rho(a)b)c = \rho(a)(bc).$$
Since $\rho(\rho(b)a) = \rho(a)b$, 
we have $$\rho(\rho(b)a)c = \rho(a)(bc).$$
Hence the restriction of $C$ to $D \setminus U$ uniquely extends to an $X$-coloring $C'$ of $D'$ 
such that the restriction of $C'$ to $D' \cap U$ corresponds to the right side of the Y-oriented R$6$(C) move depicted in Fig.~\ref{fig Y-ori'd IH-moves ver 2}. 

Conversely,
by Lemma~\ref{Lemma:Coloring condition},
for any $X$-coloring $C'$ of $D'$,
the restriction of $C'$ to $D' \cap U$ \textcolor{black}{can be assumed to be} 
as illustrated on the right side of the Y-oriented R$6$(C) move depicted in Fig.~\ref{fig Y-ori'd IH-moves ver 2},
where $a, b, c \in X$ such that
$(b, c) \in P$ and $(\rho(a), bc) \in P$.
By the associativity of the groupoid $X$, \textcolor{black}{we have}
$$(\rho(a), b) \in P\ \mbox{and}\ (\rho(a)b, c) \in P.$$
\textcolor{black}{By} Lemma~\ref{Lemma:well-def}, \textcolor{black}{remarking that $\rho$ is involutive,} 
$(\rho(a), b) \in P$ implies \textcolor{black}{that} $$(\rho(b), a) \in P\ \mbox{and}\ \rho(b)a = \rho(\rho(a)b).$$
Applying $\rho$ to both sides of $\rho(b)a = \rho(\rho(a)b)$, 
we have $$\rho(\rho(b)a) = \rho(a) b.$$
Since $(\rho(a)b, c) \in P$ and $\rho(a) b = \rho(\rho(b)a)$,
we have $$(\rho(\rho(b)a), c) \in P.$$
By the associativity of the groupoid $X$,
$$\rho(\rho(b)a)c = (\rho(a)b)c = \rho(a)(bc).$$
Hence the restriction of $C'$ to $D' \setminus U$ uniquely extends to an $X$-coloring $C$ of $D$ 
such that the restriction of $C$ to $D \cap U$ corresponds to the left side of the Y-oriented R$6$(C) move depicted in Fig.~\ref{fig Y-ori'd IH-moves ver 2}. 

Therefore,
the map $\mathrm{Col}_{X}(D) \to \mathrm{Col}_{X}(D')$ which sends $C$ to $C'$ as above is a bijection.

When $D$ and $D'$ are related by one of \textcolor{black}{the} other Y-oriented R$6$ moves depicted in Fig.~\ref{fig Y-ori'd IH-moves ver 2},
the existence of a bijection between $\mathrm{Col}_{X}(D)$ and $\mathrm{Col}_{X}(D')$ is similarly ensured by
Lemmas~\ref{Lemma:well-def} \textcolor{black}{and} \ref{Lemma:Coloring condition}, and the associativity of the groupoid $X$ (\textcolor{black}{see} Fig.~\ref{fig Y-ori'd IH-moves ver 2}).
\end{proof}

\begin{figure}[h]
  \centering
  \begin{tikzpicture}[line width=1.6pt]
    \draw (0,3)--(1,2)--(2,3);
    \draw (1,2)--(1,1);
    \draw (0,0)--(1,1)--(2,0);
    \draw (0.5,2.5) node{\rotatebox{-45}{$\blacktriangleright$}};
    \draw (1.5,2.5) node{\rotatebox{45}{$\blacktriangleleft$}};
    \draw (1,1.5) node{$\blacktriangledown$};
    \draw (0.5,0.5) node{\rotatebox{45}{$\blacktriangleright$}};
    \draw (1.5,0.5) node{\rotatebox{-45}{$\blacktriangleright$}};
    \filldraw (1,2) circle(0.15);
    \filldraw (1,1) circle(0.15);
    \draw (0,3) node[left]{$b$};
    \draw (2,3) node[right]{$c$};
    \draw (0,0) node[left]{$a$};
    \draw (1,1.5) node[left]{$bc$};
    \draw (2,0) node[right]{$a(bc)$};

    \draw[<->] (2.5,1.5)--(3.5,1.5);
    \draw (3,1) node{R$6$(A)};

    \draw (4,2.5)--(5,1.5)--(6,1.5)--(7,2.5);
    \draw (4,0.5)--(5,1.5);
    \draw (6,1.5)--(7,0.5);
    \filldraw (5,1.5) circle(0.15);
    \filldraw (6,1.5) circle(0.15);
    \draw (4.5,2) node{\rotatebox{-45}{$\blacktriangleright$}};
    \draw (4.5,1) node{\rotatebox{45}{$\blacktriangleright$}};
    \draw (5.5,1.5) node{$\blacktriangleright$};
    \draw (6.5,2) node{\rotatebox{45}{$\blacktriangleleft$}};
    \draw (6.5,1) node{\rotatebox{-45}{$\blacktriangleright$}};
    \draw (4,2.5) node[left]{$b$};
    \draw (7,2.5) node[right]{$c$};
    \draw (4,0.5) node[left]{$a$};
    \draw (5.5,1) node{$ab$};
    \draw (7,0.5) node[right]{$(ab)c$};
  \end{tikzpicture}

  \begin{tikzpicture}[line width=1.6pt]
    \draw (0,3)--(1,2)--(2,3);
    \draw (1,2)--(1,1);
    \draw (0,0)--(1,1)--(2,0);
    \draw (0.5,2.5) node{\rotatebox{-45}{$\blacktriangleright$}};
    \draw (1.5,2.5) node{\rotatebox{45}{$\blacktriangleright$}};
    \draw (1,1.5) node{$\blacktriangledown$};
    \draw (0.5,0.5) node{\rotatebox{45}{$\blacktriangleleft$}};
    \draw (1.5,0.5) node{\rotatebox{-45}{$\blacktriangleright$}};
    \filldraw (1,2) circle(0.15);
    \filldraw (1,1) circle(0.15);
    \draw (0,3) node[left]{$(ab) c$};
    \draw (2,3) node[right]{$c$};
    \draw (0,0) node[left]{$a$};
    \draw (1,1.5) node[left]{$ab$};
    \draw (2,0) node[right]{$b$};

    \draw[<->] (2.5,1.5)--(3.5,1.5);
    \draw (3,1) node{R$6$(B)};

    \draw (4,2.5)--(5,1.5)--(6,1.5)--(7,2.5);
    \draw (4,0.5)--(5,1.5);
    \draw (6,1.5)--(7,0.5);
    \filldraw (5,1.5) circle(0.15);
    \filldraw (6,1.5) circle(0.15);
    \draw (4.5,2) node{\rotatebox{-45}{$\blacktriangleright$}};
    \draw (4.5,1) node{\rotatebox{45}{$\blacktriangleleft$}};
    \draw (5.5,1.5) node{$\blacktriangleright$};
    \draw (6.5,2) node{\rotatebox{45}{$\blacktriangleright$}};
    \draw (6.5,1) node{\rotatebox{-45}{$\blacktriangleright$}};
    \draw (4,2.5) node[left]{$a(bc)$};
    \draw (7,2.5) node[right]{$c$};
    \draw (4,0.5) node[left]{$a$};
    \draw (5.5,1) node{$bc$};
    \draw (7,0.5) node[right]{$b$};
  \end{tikzpicture}

  \begin{tikzpicture}[line width=1.6pt]
    \draw (-2,2.5)--(-3,1.5)--(-4,1.5)--(-5,2.5);
    \draw (-5,0.5)--(-4,1.5);
    \draw (-2,0.5)--(-3,1.5);
    \filldraw (-3,1.5) circle(0.15);
    \filldraw (-4,1.5) circle(0.15);
    \draw (-4.5,2) node{\rotatebox{-45}{$\blacktriangleright$}};
    \draw (-4.5,1) node{\rotatebox{45}{$\blacktriangleleft$}};
    \draw (-3.5,1.5) node{$\blacktriangleleft$};
    \draw (-2.5,2) node{\rotatebox{45}{$\blacktriangleleft$}};
    \draw (-2.5,1) node{\rotatebox{-45}{$\blacktriangleright$}};
    \draw (-5,2.5) node[left]{$b$};
    \draw (-2,2.5) node[right]{$c$};
    \draw (-5,0.5) node[left]{$a$};
    \draw (-3.5,1) node{$\rho(b) a$};
    \draw (-2,0.5) node[below]{$\rho(\rho(b) a) c$};

    \draw[<->] (-0.5,1.5)--(-1.5,1.5);
    \draw (-1,1) node{R$6$(C)}; 

    \draw (0,3)--(1,2)--(2,3);
    \draw (1,2)--(1,1);
    \draw (0,0)--(1,1)--(2,0);
    \draw (0.5,2.5) node{\rotatebox{-45}{$\blacktriangleright$}};
    \draw (1.5,2.5) node{\rotatebox{45}{$\blacktriangleleft$}};
    \draw (1,1.5) node{$\blacktriangledown$};
    \draw (0.5,0.5) node{\rotatebox{45}{$\blacktriangleleft$}};
    \draw (1.5,0.5) node{\rotatebox{-45}{$\blacktriangleright$}};
    \filldraw (1,2) circle(0.15);
    \filldraw (1,1) circle(0.15);
    \draw (0,3) node[left]{$b$};
    \draw (2,3) node[right]{$c$};
    \draw (0,0) node[left]{$a$};
    \draw (1,1.5) node[left]{$bc$};
    \draw (2,0) node[right]{$\rho(a) (bc)$};

    \draw[<->] (2.5,1.5)--(3.5,1.5);
    \draw (3,1) node{R$6$(D)};

    \draw (4,2.5)--(5,1.5)--(6,1.5)--(7,2.5);
    \draw (4,0.5)--(5,1.5);
    \draw (6,1.5)--(7,0.5);
    \filldraw (5,1.5) circle(0.15);
    \filldraw (6,1.5) circle(0.15);
    \draw (4.5,2) node{\rotatebox{-45}{$\blacktriangleright$}};
    \draw (4.5,1) node{\rotatebox{45}{$\blacktriangleleft$}};
    \draw (5.5,1.5) node{$\blacktriangleright$};
    \draw (6.5,2) node{\rotatebox{45}{$\blacktriangleleft$}};
    \draw (6.5,1) node{\rotatebox{-45}{$\blacktriangleright$}};
    \draw (4,2.5) node[left]{$b$};
    \draw (7,2.5) node[right]{$c$};
    \draw (4,0.5) node[left]{$a$};
    \draw (5.5,1) node{$\rho(a) b$};
    \draw (7,0.5) node[right]{$(\rho(a) b) c$};
  \end{tikzpicture}

  \begin{tikzpicture}[line width=1.6pt]
    \draw (-2,2.5)--(-3,1.5)--(-4,1.5)--(-5,2.5);
    \draw (-5,0.5)--(-4,1.5);
    \draw (-2,0.5)--(-3,1.5);
    \filldraw (-3,1.5) circle(0.15);
    \filldraw (-4,1.5) circle(0.15);
    \draw (-4.5,2) node{\rotatebox{-45}{$\blacktriangleright$}};
    \draw (-4.5,1) node{\rotatebox{45}{$\blacktriangleleft$}};
    \draw (-3.5,1.5) node{$\blacktriangleleft$};
    \draw (-2.5,2) node{\rotatebox{45}{$\blacktriangleright$}};
    \draw (-2.5,1) node{\rotatebox{-45}{$\blacktriangleleft$}};
    \draw (-5,2.5) node[left]{$b$};
    \draw (-2,2.5) node[right]{$c$};
    \draw (-5,0.5) node[left]{$a$};
    \draw (-3.5,1) node{$\rho(b) a$};
    \draw (-2,0.5) node[below]{$c (\rho(b) a)$};

    \draw[<->] (-0.5,1.5)--(-1.5,1.5);
    \draw (-1,1) node{R$6$(E)}; 

    \draw (0,3)--(1,2)--(2,3);
    \draw (1,2)--(1,1);
    \draw (0,0)--(1,1)--(2,0);
    \draw (0.5,2.5) node{\rotatebox{-45}{$\blacktriangleright$}};
    \draw (1.5,2.5) node{\rotatebox{45}{$\blacktriangleright$}};
    \draw (1,1.5) node{$\blacktriangledown$};
    \draw (0.5,0.5) node{\rotatebox{45}{$\blacktriangleleft$}};
    \draw (1.5,0.5) node{\rotatebox{-45}{$\blacktriangleleft$}};
    \filldraw (1,2) circle(0.15);
    \filldraw (1,1) circle(0.15);
    \draw (0,3) node[left]{$b$};
    \draw (2,3) node[right]{$c$};
    \draw (0,0) node[left]{$a$};
    \draw (1,1.5) node[left]{$b\rho(c)$};
    \draw (2,0) node[right]{$\rho(b \rho(c)) a$};

    \draw[<->] (2.5,1.5)--(3.5,1.5);
    \draw (3,1) node{R$6$(F)};

    \draw (4,2.5)--(5,1.5)--(6,1.5)--(7,2.5);
    \draw (4,0.5)--(5,1.5);
    \draw (6,1.5)--(7,0.5);
    \filldraw (5,1.5) circle(0.15);
    \filldraw (6,1.5) circle(0.15);
    \draw (4.5,2) node{\rotatebox{-45}{$\blacktriangleright$}};
    \draw (4.5,1) node{\rotatebox{45}{$\blacktriangleleft$}};
    \draw (5.5,1.5) node{$\blacktriangleright$};
    \draw (6.5,2) node{\rotatebox{45}{$\blacktriangleright$}};
    \draw (6.5,1) node{\rotatebox{-45}{$\blacktriangleleft$}};
    \draw (4,2.5) node[left]{$b$};
    \draw (7,2.5) node[right]{$c$};
    \draw (4,0.5) node[left]{$a$};
    \draw (5.5,1) node{$\rho(a) b$};
    \draw (7,0.5) node[right]{$c \rho(\rho(a) b)$};
  \end{tikzpicture}
  \caption{Y-oriented R$6$ moves with $X$-colorings}{\label{fig Y-ori'd IH-moves ver 2}}
\end{figure} 

\lem{\label{Lemma:Invariance_under_inverse}}{
  Let $D$ and $D'$ be Y-oriented diagrams of spatial surfaces.
  If $D'$ is obtained from $D$ by an inverse move on an $S^{1}$-component $K$, 
  there is a bijection between $\mathrm{Col}_{X}(D)$ and $\mathrm{Col}_{X}(D')$.
}\upshape

\begin{proof}
  For any $X$-coloring $C$ of $D$,
  consider an $X$-coloring $C'$ of $D'$ defined by
  $$C'(a) = \begin{cases}
    C(a) & \mbox{\textcolor{black}{if}}\ a \ \mbox{is an arc of}\ D' \ \mbox{not lying on}\ K\textcolor{black}{,} \\
    \rho(C(a)) & \mbox{\textcolor{black}{if}}\ a \ \mbox{is an arc of}\ D' \ \mbox{lying on}\ K\textcolor{black}{.}
  \end{cases}$$
  Since $(X, \rho)$ is a symmetric rack,
  we have a bijection between $\mathrm{Col}_{X}(D)$ and $\mathrm{Col}_{X}(D')$ by sending $C$ to $C'$, \textcolor{black}{see} \cite{Kamada2007,Kamada-Oshiro2010}.
\end{proof}

\begin{proof}[Proof of Theorem~\ref{thm coloring invariant}]
  By Lemma~\ref{Lemma:generating set} and Theorem~\ref{thm Y-oriented R-moves},
  it is sufficient to consider the case that $D'$ is obtained from $D$ by a Y-oriented move depicted in Fig.~\ref{fig generating set of Y-oriented Reidemeister moves for oriented spatial surface diagrams}
  or an inverse move.

  Since $X$ is a rack,
  we see that there is a bijection between $\mathrm{Col}_{X}(D)$ and $\mathrm{Col}_{X}(D')$ in the case of a Y-oriented R$2$ move and a Y-oriented R$3$ move, \textcolor{black}{see} \cite{Fenn-Rourke1992}.
  \textcolor{black}{We saw the other cases in Lemmas~\ref{Lemma:Invariance_under_R5-moves}, \ref{Lemma:Invariance_under_R6-moves}, and \ref{Lemma:Invariance_under_inverse}.}
\end{proof}

\section{A universality of groupoid racks on colorings}{\label{Sect:Universality}}

The following theorem implies that a groupoid rack has a universal property on colorings for Y-oriented diagrams of spatial surfaces, i.e., under a certain assumption on colorings, 
any algebraic \textcolor{black}{system} that is used for colorings for Y-oriented diagrams of spatial surfaces 
has a structure of a groupoid rack.

\thm{\label{thm universality}}{
  Let $(R, \ast)$ be a rack.
  Assume that a subset $P \subset R \times R$ and a map $\mu : P \to R$
  satisfying the following conditions $(1)$--$(8)$,
  where we denote $\mu(a, b)$ by $ab$.
  \begin{enumerate}
    \item[$(1)$] For any $a, b, c \in R$, the following are equivalent.
      $$\left[(a, b) \in P \wedge (ab, c) \in P\right] \ \mbox{and} \  \left[(b, c) \in P \wedge (a, bc) \in P\right].$$
    \item[$(2)$] For any $a, b, c \in R$ with $(a, b) \in P \wedge (ab, c) \in P$, we have
      $$(ab)c = a(bc).$$
    \item[$(3)$] For any $a, b, x \in R$, the following are equivalent.
      $$(a, b) \in P \   \mbox{and} \  (a \ast x, b \ast x) \in P.$$
    \item[$(4)$] For any $(a, b) \in P$ and $x \in R$, we have
      $$(ab) \ast x = (a \ast x)(b \ast x)\ \mbox{and}\  x \ast (ab) = (x \ast a) \ast b.$$
    \item[$(5)$] For any $a, b, c, d \in R$, the following are equivalent.
      \begin{enumerate}
        \item[$(5.1)$] There exists an element $e \in R$ such that $(a, e), (e, d) \in P$, $ae = c$, and $ed = b$.
        \item[$(5.2)$] $(a, b), (c, d) \in P$ and $ab = cd$.
        \item[$(5.3)$] There exists an element $f \in R$ such that $(c, f), (f, b) \in P$, $cf = a$, and $fb = d$. 
      \end{enumerate}
    \item[$(6)$] For any $a, b, c, d \in R$, the following holds.
    \begin{enumerate}
      \item[$(6.1)$] If $(5.1)$ is true, then such an element $e$ is unique.
      \item[$(6.2)$] If $(5.3)$ is true, then such an element $f$ is unique. 
    \end{enumerate}
      
    \item[$(7)$] For any $a, b, c, d \in R$, the following are equivalent.  
      \begin{itemize}
        \item[$(7.1)$] There exists an element $x \in R$ such that $(x, b), (x, d) \in P$, $xb = a$, and $xd = c$.
        \item[$(7.2)$] There exists an element $y \in R$ such that $(a, y), (b, y) \in P$, $ay = c$, and $by = d$.
        \item[$(7.3)$] There exists an element $z \in R$ such that $(c, z), (d, z) \in P$, $cz = a$, and $dz = b$.
      \end{itemize}
    \item[$(8)$] For any $a, b, c, d \in R$, the following holds.
    \begin{enumerate}
      \item[$(8.1)$] If $(7.1)$ is true, then such an element $x$ is unique.
      \item[$(8.2)$] If $(7.2)$ is true, then such an element $y$ is unique.
      \item[$(8.3)$] If $(7.3)$ is true, then such an element $z$ is unique. 
    \end{enumerate}
  \end{enumerate}
  Put $X = \bigcup_{(a, b) \in P}\left\{a, b\right\}$ and we denote $\ast|_{X \times X}$ also by $\ast$.
  Then 
  \begin{enumerate}
    \item[$(\mathrm{i})$] $(X, \ast)$ is a subrack of the rack $(R, \ast)$.
    \item[$(\mathrm{ii})$] There is a groupoid $\mathcal{C}$ such that $(X, \ast)$ is a groupoid rack associated with $\mathcal{C}$.
  \end{enumerate}
  }

\upshape

Before giving a proof,
we show how the conditions of Theorem~\ref{thm universality} are obtained when we consider colorings of Y-oriented diagrams 
\textcolor{black}{yielding invariants of spatial surfaces.}

Let $D$ be a Y-oriented diagram and let $R = (R, \ast)$ be a pair of a set $R$ and a binary operation $\ast$ on $R$.
In this section, 
a map $\mathcal{A}(D) \to R$ is called a \textit{coloring} of $D$ by $R$ and the image of an arc $a$ by a coloring is called a \textit{color} of $a$.

Let $D$ and $D'$ be Y-oriented diagrams such that
$D'$ is obtained from $D$ by a Y-oriented Reidemeister move with support $U$.
Let $C$ and $C'$ be colorings of $D$ and $D'$ by $R$, respectively.
We say that $C$ and $C'$ \textit{are compatible} with respect to $U$
if the restriction of $C$ to $D \setminus U$ and the restriction of $C'$ to $D' \setminus U = D \setminus U$ are the same.
The \textit{coloring assumption} is the following condition:

For any coloring $C$ of $D$,
there exists a unique coloring $C'$ of $D'$ such that
$C$ and $C'$ are compatible with respect to $U$
\textcolor{black}{and for any coloring $C'$ of $D'$,
there exists a unique coloring $C$ of $D$ such that
$C$ and $C'$ are compatible with respect to $U$.}

First, we assume \textcolor{black}{that any coloring satisfies} the condition depicted in Fig.~\ref{fig coloring conditions around crossings} at \textcolor{black}{each} crossing, 
which we call the \textit{coloring condition} at \textcolor{black}{a} crossing.
\textcolor{black}{Then} we see that $R$ needs to be a rack for the coloring assumptions on Y-oriented R$2$ moves and Y-oriented R$3$ moves, \textcolor{black}{see} \cite{Fenn-Rourke1992,Joyce1982,Matveev1982}.


\begin{figure}[h]
  \centering
    \begin{tikzpicture}[line width=2pt, use Hobby shortcut]
      \begin{knot}[
        consider self intersections=true,
        ignore endpoint intersections=false,
        flip crossing/.list={}
        ]
      \strand (2,2)--(0,0);
      \strand (0,2)--(2,0);
      \end{knot}  
      \draw (0,0) node{\rotatebox{225}{$\blacktriangleright$}};
      \draw (-0.3,2) node{$x$};
      \draw (2.3,2) node{$y$};
      \draw (2.5,0) node{$x \ast y$};
    \end{tikzpicture}
    \caption{Coloring condition at a crossing ($x, y \in R$)}{\label{fig coloring conditions around crossings}}  
\end{figure}

Let $P$ be a subset of $R \times R$ and $\mu : P \to R$ a map.
We denote $\mu (a, b)$ by $ab$.
\textcolor{black}{Furthermore, assume that any coloring satisfies} the coloring conditions shown in Fig.~\ref{fig coloring condition at vertex} at \textcolor{black}{each vertex}, which we call the \textit{coloring conditions} at \textcolor{black}{a vertex}, where $(a, b) \in P$.

\begin{figure}[h]
  \centering
    \begin{tikzpicture}[line width=2pt]
      \draw (0,2)--(1,1)--(1,0);
      \draw (2,2)--(1,1);
      \draw (0.5,1.5) node{\rotatebox{315}{$\blacktriangleright$}};
      \draw (1.5,1.5) node{\rotatebox{225}{$\blacktriangleright$}};
      \draw (1,0.5) node{\rotatebox{270}{$\blacktriangleright$}};
      \draw (-0.3,2) node{$a$};
      \draw (2.3,2) node{$b$};
      \draw (1.5,0) node{$ab$};
      \filldraw (1,1) circle(0.125);
    \end{tikzpicture}
    \quad 
    \begin{tikzpicture}[line width=2pt]
      \draw (1,2)--(1,1)--(0,0);
      \draw (2,0)--(1,1);
      \draw (0.5,0.5) node{\rotatebox{225}{$\blacktriangleright$}};
      \draw (1.5,0.5) node{\rotatebox{315}{$\blacktriangleright$}};
      \draw (1,1.5) node{\rotatebox{270}{$\blacktriangleright$}};
      \filldraw (1,1) circle(0.125);
      \draw (1.5,2) node{$ab$};
      \draw (-0.3,0) node{$a$};
      \draw (2.3,0) node{$b$};
    \end{tikzpicture} 
    \caption{Coloring conditions at vertices ($(a, b) \in P$)}{\label{fig coloring condition at vertex}}  
\end{figure}

\textcolor{black}{By the assumption, we first} see that the following conditions (A.1) and (A.2) are required for the coloring assumption on Y-oriented ${\rm R}5$ moves, by observing the Y-oriented R$5$(A), R$5$(B), R$5$(C) and R$5$(D) moves \textcolor{black}{depicted in} Fig.~\ref{fig Y-oriented R$5$ moves}. 
\begin{itemize}
  \item[(A.1)] For any $a, b, x \in R$, the following are equivalent.
    $$(a, b) \in P \  \mbox{and} \  (a \ast x, b \ast x) \in P.$$
  \item[(A.2)] For any $(a, b) \in P$ and $x \in R$, we have 
    $$(ab) \ast x = (a \ast x)(b \ast x) \  \mbox{and} \  x \ast (ab) = (x \ast a) \ast b.$$
\end{itemize}

\textcolor{black}{Next}, we see that the following conditions (B.1) and (B.2) are required for the coloring assumption on Y-oriented ${\rm R}6$ moves, by observing the Y-oriented R$6$(A) and R$6$(B) moves \textcolor{black}{depicted in} Fig.~\ref{fig Y-ori'd IH-moves ver 2}.
\begin{itemize}
  \item[(B.1)] For any $a, b, c \in R$, the following are equivalent.
    $$\left[(a, b) \in P \wedge (ab, c) \in P\right] \  \mbox{and} \  \left[(b, c) \in P \wedge (a, bc) \in P\right].$$
  \item[(B.2)] For any $(a, b), (ab, c) \in P$, we have
    $$(ab)c = a(bc).$$
\end{itemize}

\textcolor{black}{Then}, we see that the following condition (C) is required for the coloring assumtion on Y-oriented ${\rm R}6$ moves, by observing the Y-oriented R$6$(C) and R$6$(D) moves \textcolor{black}{depicted in} Fig.~\ref{fig Y-ori'd IH-moves ver 3}.

\begin{itemize}
  \item[(C)] For any $a, b, c, d \in R$, the following are equivalent.
  \begin{enumerate}
    \item[(C.1)] There exists an element $e \in R$ such that $(a, e), (e, d) \in P$, $ae = c$, and $ed = b$.
    \item[(C.2)] $(a, b), (c, d) \in P$ and $ab = cd$.
    \item[(C.3)] There exists an element $f \in R$ such that $(c, f), (f, b) \in P$, $cf = a$, and $fb = d$.
  \end{enumerate}
\end{itemize} 


Furthermore for the coloring assumption,
the following conditions (D.1) and (D.2) are required.

\begin{itemize}
  \item[(D.1)] If (C.1) is true, then such an element $e$ is unique.
  \item[(D.2)] If (C.3) is true, then such an element $f$ is unique. 
\end{itemize}

\textcolor{black}{Finally,} we see that the following condition (E) is required for the coloring assumtion on Y-oriented ${\rm R}6$ moves, by observing the Y-oriented R$6$(E) and R$6$(F) moves \textcolor{black}{depicted in} Fig.~\ref{fig Y-ori'd IH-moves ver 3}.

\begin{itemize}
  \item[(E)] For any $a, b, c, d \in R$, the following are equivalent.
  \begin{enumerate}
    \item[(E.1)] There exists an element $x \in R$ such that $(x, b), (x, d) \in P$, $xb = a$, and $xd = c$.
    \item[(E.2)] There exists an element $y \in R$ such that $(a, y), (b, y) \in P$, $ay = c$, and $by = d$.
    \item[(E.3)] There exists an element $z \in R$ such that $(c, z), (d, z) \in P$, $cz = a$, and $dz = b$.
  \end{enumerate}
\end{itemize} 

Furthermore for the coloring assumption,
the following conditions (F.1), (F.2), and (F.3) are required.

\begin{itemize}
  \item[(F.1)] If (E.1) is true, then such an element $x$ is unique.
  \item[(F.2)] If (E.2) is true, then such an element $y$ is unique.
  \item[(F.3)] If (E.3) is true, then such an element $z$ is unique. 
\end{itemize}

\begin{figure}[h]
  \centering
  \begin{tikzpicture}[line width=1.6pt]
    \draw (0,3) -- (1,2) -- (2,3);
    \draw (1,2) -- (1,1);
    \draw (0,0) -- (1,1) -- (2,0);
    \draw (0.5,2.5) node{\rotatebox{-45}{$\blacktriangleright$}};
    \draw (1.5,2.5) node{\rotatebox{45}{$\blacktriangleleft$}};
    \draw (1,1.5) node{$\blacktriangledown$};
    \draw (0.5,0.5) node{\rotatebox{45}{$\blacktriangleright$}};
    \draw (1.5,0.5) node{\rotatebox{-45}{$\blacktriangleright$}};
    \filldraw (1,2) circle(0.15);
    \filldraw (1,1) circle(0.15);
    \draw (0,3) node[left]{$b$};
    \draw (2,3) node[right]{$c$};
    \draw (0,0) node[left]{$a$};
    \draw (1,1.5) node[left]{$bc$};
    \draw (2,0) node[right]{$a(bc)$};

    \draw[<->] (2.5,1.5) -- (3.5,1.5);
    \draw (3,1) node{R$6$(A)};

    \draw (4,2.5) -- (5,1.5) -- (6,1.5) -- (7,2.5);
    \draw (4,0.5) -- (5,1.5);
    \draw (6,1.5) -- (7,0.5);
    \filldraw (5,1.5) circle(0.15);
    \filldraw (6,1.5) circle(0.15);
    \draw (4.5,2) node{\rotatebox{-45}{$\blacktriangleright$}};
    \draw (4.5,1) node{\rotatebox{45}{$\blacktriangleright$}};
    \draw (5.5,1.5) node{$\blacktriangleright$};
    \draw (6.5,2) node{\rotatebox{45}{$\blacktriangleleft$}};
    \draw (6.5,1) node{\rotatebox{-45}{$\blacktriangleright$}};
    \draw (4,2.5) node[left]{$b$};
    \draw (7,2.5) node[right]{$c$};
    \draw (4,0.5) node[left]{$a$};
    \draw (5.5,1) node{$ab$};
    \draw (7,0.5) node[right]{$(ab)c$};
  \end{tikzpicture}

  \begin{tikzpicture}[line width=1.6pt]
    \draw (0,3) -- (1,2) -- (2,3);
    \draw (1,2) -- (1,1);
    \draw (0,0) -- (1,1) -- (2,0);
    \draw (0.5,2.5) node{\rotatebox{-45}{$\blacktriangleright$}};
    \draw (1.5,2.5) node{\rotatebox{45}{$\blacktriangleright$}};
    \draw (1,1.5) node{$\blacktriangledown$};
    \draw (0.5,0.5) node{\rotatebox{45}{$\blacktriangleleft$}};
    \draw (1.5,0.5) node{\rotatebox{-45}{$\blacktriangleright$}};
    \filldraw (1,2) circle(0.15);
    \filldraw (1,1) circle(0.15);
    \draw (0,3) node[left]{$(ab) c$};
    \draw (2,3) node[right]{$c$};
    \draw (0,0) node[left]{$a$};
    \draw (1,1.5) node[left]{$ab$};
    \draw (2,0) node[right]{$b$};

    \draw[<->] (2.5,1.5) -- (3.5,1.5);
    \draw (3,1) node{R$6$(B)};

    \draw (4,2.5) -- (5,1.5) -- (6,1.5) -- (7,2.5);
    \draw (4,0.5) -- (5,1.5);
    \draw (6,1.5) -- (7,0.5);
    \filldraw (5,1.5) circle(0.15);
    \filldraw (6,1.5) circle(0.15);
    \draw (4.5,2) node{\rotatebox{-45}{$\blacktriangleright$}};
    \draw (4.5,1) node{\rotatebox{45}{$\blacktriangleleft$}};
    \draw (5.5,1.5) node{$\blacktriangleright$};
    \draw (6.5,2) node{\rotatebox{45}{$\blacktriangleright$}};
    \draw (6.5,1) node{\rotatebox{-45}{$\blacktriangleright$}};
    \draw (4,2.5) node[left]{$a(bc)$};
    \draw (7,2.5) node[right]{$c$};
    \draw (4,0.5) node[left]{$a$};
    \draw (5.5,1) node{$bc$};
    \draw (7,0.5) node[right]{$b$};
  \end{tikzpicture}

  \begin{tikzpicture}[line width=1.6pt]
    \draw (-2,2.5) -- (-3,1.5) -- (-4,1.5) -- (-5,2.5);
    \draw (-5,0.5) -- (-4,1.5);
    \draw (-2,0.5) -- (-3,1.5);
    \filldraw (-3,1.5) circle(0.15);
    \filldraw (-4,1.5) circle(0.15);
    \draw (-4.5,2) node{\rotatebox{-45}{$\blacktriangleright$}};
    \draw (-4.5,1) node{\rotatebox{45}{$\blacktriangleleft$}};
    \draw (-3.5,1.5) node{$\blacktriangleleft$};
    \draw (-2.5,2) node{\rotatebox{45}{$\blacktriangleleft$}};
    \draw (-2.5,1) node{\rotatebox{-45}{$\blacktriangleright$}};
    \draw (-5,2.5) node[left]{$a$};
    \draw (-2,2.5) node[right]{$b$};
    \draw (-5,0.5) node[left]{$c$};
    \draw (-3.5,1) node{$e$};
    \draw (-2,0.5) node[below]{$d$};

    \draw[<->] (-0.5,1.5) -- (-1.5,1.5);
    \draw (-1,1) node{R$6$(C)}; 

    \draw (0,3) -- (1,2) -- (2,3);
    \draw (1,2) -- (1,1);
    \draw (0,0) -- (1,1) -- (2,0);
    \draw (0.5,2.5) node{\rotatebox{-45}{$\blacktriangleright$}};
    \draw (1.5,2.5) node{\rotatebox{45}{$\blacktriangleleft$}};
    \draw (1,1.5) node{$\blacktriangledown$};
    \draw (0.5,0.5) node{\rotatebox{45}{$\blacktriangleleft$}};
    \draw (1.5,0.5) node{\rotatebox{-45}{$\blacktriangleright$}};
    \filldraw (1,2) circle(0.15);
    \filldraw (1,1) circle(0.15);
    \draw (0,3) node[left]{$a$};
    \draw (2,3) node[right]{$b$};
    \draw (0,0) node[left]{$c$};
    \draw (1,1.35) node[left]{$ab = cd$};
    \draw (2,0) node[right]{$d$};

    \draw[<->] (2.5,1.5) -- (3.5,1.5);
    \draw (3,1) node{R$6$(D)};

    \draw (4,2.5) -- (5,1.5) -- (6,1.5) -- (7,2.5);
    \draw (4,0.5) -- (5,1.5);
    \draw (6,1.5) -- (7,0.5);
    \filldraw (5,1.5) circle(0.15);
    \filldraw (6,1.5) circle(0.15);
    \draw (4.5,2) node{\rotatebox{-45}{$\blacktriangleright$}};
    \draw (4.5,1) node{\rotatebox{45}{$\blacktriangleleft$}};
    \draw (5.5,1.5) node{$\blacktriangleright$};
    \draw (6.5,2) node{\rotatebox{45}{$\blacktriangleleft$}};
    \draw (6.5,1) node{\rotatebox{-45}{$\blacktriangleright$}};
    \draw (4,2.5) node[left]{$a$};
    \draw (7,2.5) node[right]{$b$};
    \draw (4,0.5) node[left]{$c$};
    \draw (5.5,1) node{$f$};
    \draw (7,0.5) node[right]{$d$};
  \end{tikzpicture}

  \begin{tikzpicture}[line width=1.6pt]
    \draw (-2,2.5) -- (-3,1.5) -- (-4,1.5) -- (-5,2.5);
    \draw (-5,0.5) -- (-4,1.5);
    \draw (-2,0.5) -- (-3,1.5);
    \filldraw (-3,1.5) circle(0.15);
    \filldraw (-4,1.5) circle(0.15);
    \draw (-4.5,2) node{\rotatebox{-45}{$\blacktriangleright$}};
    \draw (-4.5,1) node{\rotatebox{45}{$\blacktriangleleft$}};
    \draw (-3.5,1.5) node{$\blacktriangleleft$};
    \draw (-2.5,2) node{\rotatebox{45}{$\blacktriangleright$}};
    \draw (-2.5,1) node{\rotatebox{-45}{$\blacktriangleleft$}};
    \draw (-5,2.5) node[left]{$a$};
    \draw (-2,2.5) node[right]{$b$};
    \draw (-5,0.5) node[left]{$c$};
    \draw (-3.5,1) node{$y$};
    \draw (-2,0.5) node[below]{$d$};

    \draw[<->] (-0.5,1.5) -- (-1.5,1.5);
    \draw (-1,1) node{R$6$(E)}; 

    \draw (0,3) -- (1,2) -- (2,3);
    \draw (1,2) -- (1,1);
    \draw (0,0) -- (1,1) -- (2,0);
    \draw (0.5,2.5) node{\rotatebox{-45}{$\blacktriangleright$}};
    \draw (1.5,2.5) node{\rotatebox{45}{$\blacktriangleright$}};
    \draw (1,1.5) node{$\blacktriangledown$};
    \draw (0.5,0.5) node{\rotatebox{45}{$\blacktriangleleft$}};
    \draw (1.5,0.5) node{\rotatebox{-45}{$\blacktriangleleft$}};
    \filldraw (1,2) circle(0.15);
    \filldraw (1,1) circle(0.15);
    \draw (0,3) node[left]{$a$};
    \draw (2,3) node[right]{$b$};
    \draw (0,0) node[left]{$c$};
    \draw (1,1.5) node[left]{$x$};
    \draw (2,0) node[right]{$d$};

    \draw[<->] (2.5,1.5) -- (3.5,1.5);
    \draw (3,1) node{R$6$(F)};

    \draw (4,2.5) -- (5,1.5) -- (6,1.5) -- (7,2.5);
    \draw (4,0.5) -- (5,1.5);
    \draw (6,1.5) -- (7,0.5);
    \filldraw (5,1.5) circle(0.15);
    \filldraw (6,1.5) circle(0.15);
    \draw (4.5,2) node{\rotatebox{-45}{$\blacktriangleright$}};
    \draw (4.5,1) node{\rotatebox{45}{$\blacktriangleleft$}};
    \draw (5.5,1.5) node{$\blacktriangleright$};
    \draw (6.5,2) node{\rotatebox{45}{$\blacktriangleright$}};
    \draw (6.5,1) node{\rotatebox{-45}{$\blacktriangleleft$}};
    \draw (4,2.5) node[left]{$a$};
    \draw (7,2.5) node[right]{$b$};
    \draw (4,0.5) node[left]{$c$};
    \draw (5.5,1) node{$z$};
    \draw (7,0.5) node[right]{$d$};
  \end{tikzpicture}
  \caption{Y-oriented R$6$ moves with $R$-colorings}{\label{fig Y-ori'd IH-moves ver 3}}
\end{figure}

From the above,
we obtain the conditions $(1)$--$(8)$ of Theorem \ref{thm universality}.
Therefore, whenever we consider colorings using the coloring conditions shown in Figs. \ref{fig coloring conditions around crossings} 
and \ref{fig coloring condition at vertex}, the conditions $(1)$--$(8)$ of Theorem \ref{thm universality} are required.

In what follows, 
suppose that $R = (R, \ast)$ \textcolor{black}{is a rack and a pair of} $P \subset R \times R$ and $\mu : P \to R$ \textcolor{black}{satisfies} the conditions $(1)$--$(8)$ of Theorem \ref{thm universality}. 
Let $X = \bigcup_{(a, b) \in P}\left\{a, b\right\}$.





\lem{\label{Lemma:identity}}{
    \textcolor{black}{For any $a \in X$, there exist unique elements $1_{s_{a}}, 1_{t_{a}} \in R$ such that $$(1_{s_{a}}, a),\  (a, 1_{t_{a}}) \in P,\ 1_{s_{a}} a = a,\ \mbox{and}\ a 1_{t_{a}} = a.$$}
}\upshape

\begin{proof}
  We first show the existence of $1_{s_{a}}$ and $1_{t_{a}}$.
  By the definition of $X$, there exists an element $b \in R$ such that 
  $(a, b) \in P$ or $(b, a) \in P$.

  Assume \textcolor{black}{that} $(a, b) \in P$.
  By the condition $(5.2) \Rightarrow (5.1)$ with $(a, b, c, d) = (a, b, a, b)$,
  $(a, b) \in P$ and $ab = ab$ imply that
  there exists an element $1_{t_{a}} \in R$ such that 
  $$(a, 1_{t_{a}}) \in P \ \mbox{and}\   a1_{t_{a}} = a.$$
  



  \textcolor{black}{Furthermore,} by the condition $(7.2) \Rightarrow (7.1)$ with $(a, b, c, d, y) = (a, a, a, a, 1_{t_{a}})$,
  $(a, 1_{t_{a}}) \in P$ and $a 1_{t_{a}} = a$ imply that there exists an element $1_{s_{a}} \in R$ such that
  $$(1_{s_{a}}, a) \in P \ \mbox{and} \ 1_{s_{a}} a = a.$$



  Assume \textcolor{black}{that} $(b, a) \in P$.
  By the condition $(5.2) \Rightarrow (5.3)$ with $(a, b, c, d) = (b, a, b, a)$,
  $(b, a) \in P$ and $ba = ba$ imply that
  there exists an element $1_{s_{a}} \in R$ such that
  $$(1_{s_{a}}, a) \in P \ \mbox{and}\  1_{s_{a}} a = a.$$
  

  \textcolor{black}{Furthermore,} by the condition $(7.1) \Rightarrow (7.2)$ with $(a, b, c, d, x) = (a, a, a, a, 1_{s_{a}})$,
  $(1_{s_{a}}, a) \in P$ and $1_{s_{a}} a = a$ imply that
  there exists an element $1_{t_{a}} \in R$ such that
  $$(a, 1_{t_{a}}) \in P\ \mbox{and}\ a 1_{t_{a}} = a.$$

  Now we show the uniqueness of $1_{s_{a}}$ and $1_{t_{a}}$.
  Since we have $(1_{s_{a}}, a) \in P$ and $1_{s_{a}} a = a$,
  the condition $(7.1)$ with $(a, b, c, d, x) = (a, a, a, a, 1_{s_{a}})$ is satisfied.
  By $(8.1)$,
  such an element $1_{s_{a}}$ is unique.

  \textcolor{black}{Furthermore,} since we have $(a, 1_{t_{a}}) \in P$ and $a 1_{t_{a}} = a$,
  the condition $(7.2)$ with $(a, b, c, d, y) = (a, a, a, a, 1_{t_{a}})$ is satisfied.
  By $(8.2)$,
  such an element $1_{t_{a}}$ is unique.
\end{proof}

\lem{\label{Lemma:Inverse}}{
  For any $a \in X$, 
    there exists a unique element $r_{a} \in R$ such that
    $$(a, r_{a}), (r_{a}, a) \in P,\ ar_{a} = 1_{s_{a}}, \mbox{and}\ r_{a}a = 1_{t_{a}}.$$
}\upshape

\begin{proof}
  \textcolor{black}{By} Lemma~\ref{Lemma:identity}, we have
  $$(a, 1_{t_{a}}),\  (1_{s_{a}}, a) \in P,\ \mbox{and}\ a {1_{t_{a}}} = 1_{s_{a}} a \ ( = a).$$
  By the conditions $(5.2) \Rightarrow (5.1)$ with $(a, b, c, d) = (a, 1_{t_{a}}, 1_{s_{a}}, a)$ and $(6.1)$,
  $(a, 1_{t_{a}}), (1_{s_{a}}, a) \in P$ and $a {1_{t_{a}}} = 1_{s_{a}} a$ imply that
  there exists a unique element $r_{a} \in R$ such that
  \begin{equation*}
    (a, r_{a}),\  (r_{a}, a) \in P,\ ar_{a} = 1_{s_{a}},\ \mbox{and}\ r_{a}a = 1_{t_{a}}. \qedhere
  \end{equation*}

  
  
\end{proof}

For each $a \in X$, 
we \textcolor{black}{let $a^{-1}$ denote $r_{a}$.}

\lem{\label{Lemma:Well-definedness}}{
  \textcolor{black}{For any $(a, b) \in P$, it holds that $1_{t_{a}} = 1_{s_{b}}$.}
}
\begin{proof}
  \textcolor{black}{For any $(a, b) \in P$,
  by $(5)$ and $(6)$ with $(a, b, c, d) = (a, b, a, b)$,
  there exists an unique element $e \in R$ such that $ae = a$ and $eb = b$.
  On the other hand,
  by Lemma~\ref{Lemma:identity},
  $ae = a$ implies that $1_{t_{a}} = e$.
  Moreover,
  by Lemma~\ref{Lemma:identity},
  $eb = b$ implies that $1_{s_{b}} = e$.
  Therefore,
  $1_{t_{a}} = 1_{s_{b}}$.}
\end{proof}

\begin{proof}[Proof of Theorem \ref{thm universality}]
  (i) We show that $X$ is closed under the operations $\ast$ and $\ast^{-1}$.
  By the condition $(3)$, 
  for any $x \in X$ and $y \in R$,
  $x \ast y \in X$.
  Hence, $X$ is closed under $\ast$.
  By the condition $(3)$,
  for any $x \in X$ and $y \in R$,
  $(x \ast^{-1} y) \ast y = x$ implies that $x \ast^{-1} y \in X$.
  Hence, $X$ is closed under $\ast^{-1}$.
  Therefore, $(X, \ast)$ is a subrack of $(R, \ast)$.

  (ii) For each element $x \in X$,
  let us introduce two symbols $s_{x}$ and $t_{x}$.
  Let $W$ be the set $\bigcup_{(a, b) \in P} \left\{s_{a}, s_{b}, t_{a}, t_{b}\right\}$
  and $\sim$ 
  the equivalence relation on $W$ generated by
  $\left\{(t_{a}, s_{b}) \mid (a, b) \in P\right\} \subset W \times W$.
  We consider the groupoid $\mathcal{C}$ defined by the following.
  \begin{itemize}
    \item $\textrm{Ob}(\mathcal{C}) = W/\sim$.
    \item For any $x, y \in \textrm{Ob}(\mathcal{C})$, $\textrm{Hom}(x, y) = \left\{a \in X \mid \textcolor{black}{\left[s_{a}\right] = x, \left[t_{a}\right] = y}\right\}$.
    \item The composition $\mathrm{Hom}(x, y) \times \mathrm{Hom}(y, z) \to \mathrm{Hom}(x, z)$ is defined by $(a, b) \mapsto ab$.
    \item Let $x \in \mathrm{Ob}(\mathcal{C})$. 
    If $x = \left[s_{a}\right]$ then
    \textcolor{black}{we identify $1_{s_{a}}$ with} the identity for $x \in \textrm{Ob}(\mathcal{C})$.
    If $x = \left[t_{a}\right]$ then
    \textcolor{black}{we identify $1_{t_{a}}$ with} the identity for $x \in \textrm{Ob}(\mathcal{C})$.
    \item For any $a \in \textrm{Hom}(x, y)$, $a^{-1} (= r_{a})$ is the inverse of $a$.
  \end{itemize}

  The associativity of the composition of morphisms of $\mathcal{C}$ is ensured by the conditions $(1)$ and $(2)$.
  \textcolor{black}{The identity morphisms are well-defined by Lemma~\ref{Lemma:Well-definedness}. }
  



  Therefore, $X$ as a set \textcolor{black}{coincides with} the set of all morphisms of $\mathcal{C}$ and  
  $\mu$ is regarded as the composition of morphisms.  

  We now verify that $(X, \ast)$ satisfies the conditions (i)--(iii) of Defnition~\ref{def:groupoindrack}.

  (i) For any $x \in X$ and $a, b \in X$ with $(a, b) \in P$,
  by the condition $(4)$,
  we have $$x \ast (ab) = (x \ast a) \ast b.$$
  
  For any $a, x \in X$,
  by the condition \textcolor{black}{$(4)$ and Lemma~\ref{Lemma:identity}}, 
  \textcolor{black}{we have $$(x \ast 1_{s_{a}}) \ast a = x \ast (1_{s_{a}} a) = x \ast a.$$}
  Since $S_{a}$ is injective,
  $(x \ast 1_{s_{a}}) \ast a = x \ast a$ implies that
  $$x \ast 1_{s_{a}} = x.$$

  From \textcolor{black}{Lemmas~\ref{Lemma:Inverse} and \ref{Lemma:Well-definedness}},
  $(a, a^{-1}) \in P$ implies that $$1_{s_{a^{-1}}} = 1_{t_{a}}.$$
  Therefore, from the above argument,
  we have $$x \ast 1_{t_{a}} = x \ast 1_{s_{a^{-1}}} = x.$$

  (ii) Since $\ast = \ast|_{X \times X} : X \times X \to X$ is the restriction of the rack operation $\ast : R \times R \to R$,
  for any $x, y, z \in X$, we have 
  $$(x \ast y) \ast z = (x \ast z) \ast (y \ast z).$$

  (iii) For any $x \in X$ and $a, b \in X$ with $(a, b) \in P$,
  by the conditions $(3)$ and $(4)$,
  we have $$(a \ast x, b \ast x) \in P\ \mbox{and}\ (ab) \ast x = (a \ast x) (b \ast x).$$

  Therefore $(X, \ast)$ is a groupoid rack associated with $\mathcal{C}$.
\end{proof}

\rem{
  When we use an algebraic \textcolor{black}{system} to consider 
  a coloring of a diagram of a spatial surface with the coloring conditions of Section~\ref{Sect:Universality}, 
  the algebraic \textcolor{black}{system} needs to have a groupoid rack structure for the number of colorings to give an invariant of the spatial surface.
}

\section*{Acknowledgement}

The author would like to thank Seiichi Kamada, Atsushi Ishii, and Yuta Taniguchi for \textcolor{black}{their} helpful advice and discussions on this research. 
\textcolor{black}{He also sincerely appreciates the referee's careful reading and constructive comments, which have helped improve the paper.}




\bibliographystyle{jplain}
\bibliography{grpd.bib}
\vspace{-3mm}

\address{(K. Arai) Department of Mathematics, Graduate School of Science, Osaka University, 1-1, Machikaneyama, Toyonaka, Osaka, 560-0043, Japan}

\email{u068111h@ecs.osaka-u.ac.jp}


\end{document}